\documentclass[10pt]{amsart}
\usepackage{amsfonts}
\usepackage{color}
\usepackage[square,compress,comma, sort, numbers]{natbib}
\usepackage[colorlinks=true, citecolor=blue, linkcolor=blue]{hyperref}
\allowdisplaybreaks[4]
\usepackage{amssymb}
\usepackage{amsmath}

\definecolor{c20}{rgb}{0.,0.7,0.}
\definecolor{c30}{rgb}{0,0,0}
\definecolor{c40}{rgb}{1,0.1,0.7}
\definecolor{c50}{rgb}{0,0,0}
\definecolor{c60}{rgb}{1,0.9,0.1}

\allowdisplaybreaks[4]

\def\ehe#1{\textcolor{c50}{#1}}

\def\eE#1{\textcolor{c50}{#1}}

\def\cLa#1{\textcolor{c50}{#1}}
\def\cLb#1{\textcolor{c30}{#1}}
\def\phu{p_h(u)}

\newcommand{\kb}[1]{\boldsymbol{#1}}
\newcommand{\vk}[1]{\kb{#1}}

\newcommand{\ve}{\varepsilon}

\newcommand{\abs}[1]{\left\lvert #1 \right\rvert}

\newcommand{\norm}[1]{\lVert #1 \rVert}

\newcommand{\E}[1]{\mathbb{E}\left\{#1\right\}}

\newcommand{\pk}[1]{\mathbb{P} \left \{#1 \right \} }

\newcommand{\R}{\mathbb{R}}

\newcommand{\N}{\mathbb{N}}
\newcommand{\inr}{\in \R}

\newcommand{\ldot}{,\ldots,}

\newcommand{\limit}[1]{\lim_{#1 \to   \infty}}

\newcommand{\BQN}{\begin{eqnarray}}
\newcommand{\EQN}{\end{eqnarray}}
\newcommand{\BQNY}{\begin{eqnarray*}}
\newcommand{\EQNY}{\end{eqnarray*}}

\newcommand{\BS}{\begin{sat}}
\newcommand{\ES}{\end{sat}}
\newcommand{\BT}{\begin{theo}}
\newcommand{\ET}{\end{theo}}
\newcommand{\BL}{\begin{lem}}
\newcommand{\EL}{\end{lem}}
\newcommand{\BK}{\begin{korr}}
\newcommand{\EK}{\end{korr}}

\newcommand{\BD}{\begin{de}}
\newcommand{\ED}{\end{de}}

\newcommand{\BIT}{\begin{itemize}}
\newcommand{\EIT}{\end{itemize}}
\newcommand{\BDI}{\begin{description}}
\newcommand{\EDI}{\end{description}}

\newcommand{\BRM}{\begin{remarks}}
\newcommand{\ERM}{\end{remarks}}

\newcommand{\BEL}{\begin{lem}}
\newcommand{\EEL}{\end{lem}}

\newtheorem{theo}{Theorem}[section]
\newtheorem{sat}[theo]{Proposition}
\newtheorem{de}[theo]{Definition}
\newtheorem{lem}[theo]{Lemma}
\newtheorem{exxa}[theo]{Example}

\newtheorem{korr}[theo]{Corollary}
\newtheorem{remark}[theo]{Remark}
\newtheorem{remarks}[theo]{Remarks}

\newtheorem{cond}[theo]{Condition}

\newcommand{\nelem}[1]{{Lemma \ref{#1}}}

\newcommand{\netheo}[1]{{Theorem \ref{#1}}}

\newcommand{\necon}[1]{{Condition \ref{#1}}}

\newcommand{\prooftheo}[1]{ \textsc{\bf Proof of Theorem} \ref{#1}:}

\newcommand{\prooflem}[1]{\textsc{\bf Proof of Lemma} \ref{#1}:}
\newcommand{\proofkorr}[1]{\textsc{\bf Proof of Corollary} \ref{#1}:}

\newcommand{\COM}[1]{}

\newcommand{\QED}{\hfill $\Box$}

\def\rw{\rightarrow}

\def\IF{\infty}

\def\LT{\left}
\def\RT{\right}

\def\rw{\rightarrow}

\def\vn{\varepsilon}
\def\Var{\text{Var}}

\def\Bu+#1{\mathcal{B}^{\varepsilon+}_{u}(#1)}

\def\uy{u_{k,y}}
\def\pu#1{\cLb{\mathcal{P}_u\LT(#1\RT)}}
\def\qu#1{\cLb{\mathcal{Q}_u\LT(#1\RT)}}

\begin{document}

\title{Extremes of Gaussian chaos processes with Trend}

\author{Long Bai}
\address{Long Bai,
Department of Actuarial Science, University of Lausanne, UNIL-Dorigny, 1015 Lausanne, Switzerland
}
\email{Long.Bai@unil.ch}

\bigskip

 \maketitle

{\bf Abstract:} Let $\boldsymbol{X}(t)=(X_1(t),\ldots,X_d(t)), t\in [0,S]$ be a Gaussian vector process and let $g(\boldsymbol{x}),\boldsymbol{x}\in\mathbb{R}^d$ be a continuous homogeneous function. In this paper we are concerned with the exact tail asymptotics of the chaos process $g(\boldsymbol{X}(t))+ h(t),t\in [0,S]$ with trend function $h$.  Both scenarios $\boldsymbol{X}(t)$ is locally-stationary and  $\boldsymbol{X}(t)$ is non-stationary are considered. Important examples include the product of Gaussian processes and chi-processes.

{\bf Keywords:}  Gaussian chaos; Gaussian vector processes; asymptotic methods; Pickands constant.

{\bf AMS Classification:} Primary 60G15; secondary 60G70

\def\gg{\eE{\hat{g}}}

\section{Introduction}
Let $\vk{X}(t)=(X_1(t)\ldots,X_d(t)), t\in[0,S],d\ge 1$ be an ${\R}^d$-valued  \ehe{centered Gaussian vector process with continuous sample paths} and let $g(\vk{x}),\vk{x}\in\R^d$ be a homogeneous function of order $p>0$, i.e., $g(c \vk x)= c^p g(\vk x), p  \in (0,\IF) $
for any $c>0, \vk x \inr^d$. Adopting the terminology of \cite{PiterChaos2015}, we shall refer  to $Y(t),t\in [0,S]$ i as the Gaussian chaos process of $\vk{X}$ with respect to $g$ \ehe{defined by}
$$ Y(t)=g( \vk{X}(t)),\quad  t\in [0,S].$$
\ehe{Throughout in the following we shall assume that $g$ is not-negative}, i.e., $g(\vk x_0)>0$ for some $\vk x_0\inr^d$.\\
The tail asymptotics of the supremum of Gaussian chaos processes has been recently derived in \cite{PiterChaos2015, ZhdPit18}, see also \cite{Zhdanov} for a simpler case. Specifically, if
 $X_i$'s are centered stationary and independent Gaussian processes with unit variance and common correlation function $r$ satisfying
\BQN\label{rr}
r(t)=1-a|t|^{\alpha} +o(|t|^{\alpha}),\ \ t\rw 0,\ \  a>0,\ \alpha\in(0,2],\ \  r(t)<1,\ \  t\in(0,T],
\EQN
then by \cite{PiterChaos2015}  (under some restrictions on $g$)
\BQN \label{e1}
\pk{ \sup_{t\in [0,T]} Y(t)> u} \sim   \mathcal{H}_\alpha a^{1/\alpha} T (u/\gg)^{\frac{2}{\alpha p}}\pk{ Y(0)> u}, \quad u\to \IF ,
\EQN
where $\gg=\max_{\vk{e}\in\mathbb{S}_{d-1}}g(\vk{e})$ with $\mathbb{S}_{d-1}$  the unit sphere on $\R^d$ and \ehe{ $\mathcal{H}_\alpha=:\mathcal{H}_\alpha^0 $ the Pickands constant defined for any $\alpha, \delta$ positive by}
$$\mathcal{H}_{\alpha}^\delta =\lim_{T\rightarrow\IF}\frac{1}{T} 
\ehe{\E{\sup_{t\in \delta Z  \cap [0,T] }
	e^{\sqrt{2}B_\alpha(t)-|t|^\alpha}}},
$$
where $B_\alpha(t),t\in\mathbb{R}$ is a standard fractional Brownian motion (fBm)  with Hurst index $\alpha/2\in(0,1]$ and we interpret
$\delta Z $ as $\R$ if $\delta =0$.

In this paper, we are interested in the tail asymptotics of supremum of the Gaussian chaos process $Y(t)=g(\vk{X}(t))$ where $X_i$'s are some general non-stationary Gaussian processes considering further a trend function $h$, i.e.,  we shall investigate
\BQNY
\phu :=
\pk{\sup_{t\in[0,T]}( Y(t) +h(t))>u}
\EQNY
as $u\rw\IF$. The non-stationary case treated here is quite different from the stationary one already dealt with in the aforementioned reference. Since we allow for a trend function $h$, our results are new even for $X_i$'s being stationary.\\
The main challenges when dealing with the Gaussian chaos process $Y$ is that the Slepian inequality (see \cite{Pit96} and \cite{GennaSlepian}) does not hold in general. In the particular case of chi-square processes, using the duality of the norms, the problem in question can be related to that of supremum of a Gaussian random field (see e.g., \cite{Pit96,ELPNT2017}) and thus Slepian inequality or some modifications of it (see \cite{AdlerTaylor,Tabis}) can still be used. \\
Organisation of the rest of the paper: In Section 2 we show our main results and some examples. Following are the proofs and some useful lemmas in Section 3 and Section 4, respectively.

\section{Main results}
Next, we introduce some restriction on $g$, assuming first that
\BQN\label{gne}
\{\vk{x}:g(\vk{x})>0, \vk{x}\in \R^d\}\neq \emptyset.
\EQN
Consider spherical coordinates $\vk{x}=(r,\vk{\varphi}),$ with $\vk{\varphi}=(\varphi_1,\ldots,\varphi_{d-1})\in \Pi_{d-1}=[0,\pi)^{d-2}\times[0,2\pi),$ being the angular coordinates of $\vk{v}$, and with $r=|\vk{v}|$.
The Jacobian of the mapping is given by
\BQN\label{JJ1}
J(r,\vk{\varphi})=r^{d-1}\sin^{d-2}\varphi_1\ldots \sin\varphi_{d-2}.
\EQN
Write below $g(\vk{\varphi})=g(\vk{x}/|\vk{x}|)$ and without loss of generalities set
$$\gg:=\max_{\vk{\varphi}\in\Pi_{d-1}}g(\vk{\varphi})
=\max_{\vk{e}\in\mathbb{S}_{d-1}}g(\vk{e}),$$
$$\mathcal{M}:=\{\vk{e}\in\mathbb{S}_{d-1}:g(\vk{e})=1\},\quad \mathcal{M}_{\vk{\varphi}}:=\{\vk{\varphi}\in\Pi_{d-1}:g(\vk{\varphi})=1\}.$$
In the following we shall assume the following condition:
\begin{cond}\label{gcon}
	There exists an $\vn>0$ such that  $g(\vk{\varphi})$ is twice continuously differentiable in the neighbourhood
	\BQNY
	\mathcal{M}_{\vk{\varphi}}(\vn):=\{\vk{\varphi}\in\Pi_{d-1}:g(\vk{\varphi})>\gg-\vn\}
	\supseteq\mathcal{M}_{\vk{\varphi}}
	\EQNY
	Further we consider two different structures of the $\mathcal{M}$:\\
	(i) $\mathcal{M}$ consists of a finite number of points, $m=0$, and $|\det g''(\vk{\varphi})|>0$ for every $\vk{\varphi}\in\mathcal{M}_{\vk{\varphi}}\eE{(\ve)}$, where
	\BQN \label{g2}
	g''(\vk{\varphi}):=\LT[\frac{\partial^2g(\vk{\varphi})}{\partial \varphi_i\partial \varphi_l}\RT]_{i,l=1,\ldots,d-1}
	\EQN
	is the Hessian matrix of $g(\vk{\varphi})$.\\
	(ii) $\mathcal{M}$ is a smooth m-dimensional manifold, $1\leq m\leq d-2$, and the rank of the matrix $g''(\vk{\varphi})$ equals $d-1-m$ for all $\vk{\varphi}\in\mathcal{M}_{\vk{\varphi}}(\ve)$.
\end{cond}
Next through this paper, we always assume that $g$ satisfied \necon{gcon}.\\
\cLa{Since $h(t), t\in[0,T]$ is continuous and for $\gg\neq 1$ and $u$ large
\BQNY
\pk{\sup_{t\in[0,T]}( Y(t) +h(t))>u}=\pk{\sup_{t\in[0,T]}\LT( \LT(\frac{1}{\gg}g( \vk{X}(t))\vee 0\RT) +\frac{1}{\gg}h(t)\RT)>\frac{u}{\gg}},
\EQNY
and $g_1(\vk{x}):=\frac{1}{\gg}g(\vk{x})\vee 0$ with $\max_{\vk{e}\in\mathbb{S}_{d-1}}g_1(\vk{e})=1$ satisfies \necon{gcon} if $g$ satisfies it.\\
Then next without loss of generalities, \ehe{hereafter  we shall assume that}  $g$ is non-negative  and $\gg=1$.}\\

\subsection{Non-stationary cases}
In this section, we assume that $X_i(t), i=1,\ldots, d$, are i.i.d. centered Gaussian processes with variance function $\sigma^2(t)$ and correlation function $r(s,t)$. Further, we assume that $\sigma(t)$ attains its maximum at point $t_0\in[0,T]$ over $[0,T]$ with
\begin{align}\label{cc2}
\sigma(t)=1-b|t-t_0|^\beta, \quad t\rw t_0, b,\ \beta>0,
\end{align}
and
\BQN\label{lrr2}
r(t,s)=1-a|t-s|^{\alpha} +o(|t-s|^{\alpha}),\ \ t,s\rw t_0,\ a>0,\ \alpha\in(0,2].
\EQN
Moreover, for a continuous function $ h(t), t\in[0,T]$, we assume that
\BQN \label{ht2}
	h(t_0) - h(t) \sim c|t-t_0|^\gamma,\ \ \ \ \ t \rw t_0,
	\EQN
where $c\geq 0$ and $\gamma>0$ are some positive constants. \\
\ehe{In order to avoid unnecessary notation we shall assume in the following that $t_0=0$ or $t_0 =T$. In the non-stationary case, another important constant appearing in the asymptotics of supremum of Gaussian processes is the Piterbarg
	constant defined for some set $E\subset \R$ and $a>0, \alpha \in (0,2], \delta >0$ by
	\BQNY
	\mathcal{P}_{\alpha,a}^{f}  E =
	\E{  \sup_{t\in  E }
		e^{ \sqrt{2a}B_\alpha(t)-a|t|^\alpha-f(t)}}
	\EQNY
	see \cite{PicandsA, Pit72, debicki2002ruin,DI2005,DE2014,DiekerY, DM, SBK, GeneralPit16} for various properties of $\mathcal{H}_{\alpha}$ and $\mathcal{P}_{\alpha,a}^f$.
}

\COM{
	Define below a finite constant appearing in our results as follows
\BQNY
\widetilde{\mathcal{P}}_{\alpha,a}^{f}=\LT\{
\begin{array}{ll}
\mathcal{P}_{\alpha,a}^{f}[0,\IF), & \hbox{if} \ t_0=0\  \hbox{and}\ t_0=T,\\
\mathcal{P}_{\alpha,a}^{f}(-\IF,\IF), & \hbox{if} \ t_0\in(0,T).
\end{array}
\RT.
\EQNY
}
\BT\label{Thm3}
Suppose that $\vk{X}$ satisfies \eqref{cc2} and \eqref{lrr2} and the continuous function $h$ satisfies \eqref{ht2}. Further assume that $g$ is homogeneous of order $p\in(0,\IF)$ satisfying \eqref{gne} and set $\alpha^*=\alpha p$ and $ \beta^*=\min(\beta p,\frac{2\gamma p}{2-p})\mathbb{I}_{\{p<2\}}+\beta p\mathbb{I}_{\{p\geq2\}}$. Then we have
\BQNY
\phu \sim \mathbf{C}_{t_0}u^{(\frac{2}{\alpha^*}-\frac{2}{\beta^*})_{+}}\pk{Y(t_0)>u-h(t_0)},
\EQNY
where
\BQN\label{Ct0}
\mathbf{C}_{t_0}=\LT\{
\begin{array}{ll}\mathcal{H}_{\alpha}
a^{\frac{1}{\alpha}}\int_{0}^\IF e^{-f(t)}dt,& \text{if}\ \beta^*>\alpha^*,\\
\ehe{{\mathcal{P}}_{\alpha,a}^{f}[0,\IF) }, & \text{if}\ \beta^*=\alpha^*,\\
1,& \text{if}\ \beta^*<\alpha^*,
\end{array}
\RT.
\EQN
with $f(t)=\frac{c}{p}t^\gamma\mathbb{I}_{\{\beta^*=\frac{2\gamma p}{2-p}\}} +b t^\beta\mathbb{I}_{\{\beta^*=\beta p\}} , t\ge 0.$
\ET
\begin{remarks}
i) In \netheo{Thm3}, \ehe{as can be seen by its proof, } when $p\in[2,\IF)$ the results still hold without  \eqref{ht2}. Further, we notice that
\BQNY
\pk{Y(t_0)>u-h(t_0)}\sim\pk{Y(t_0)>u}
\LT\{
\begin{array}{ll}
e^{\frac{h(t_0)}{2}},& \text{if}\ p=2,\\
1,& \text{if}\ p>2.
\end{array}
\RT.
\EQNY
\eE{ii) In \netheo{Thm3}, setting $h(t_0)=0$ and $c=0$, we retrieve the result for  $h(t)=0, t\in [0,T]$.}\\
\ehe{iii) If $t_0 \in (0,T)$, then the above results hold with $2 \mathcal{H} _\alpha$ instead of $\mathcal{H}_\alpha$ and
	 $\ehe{{\mathcal{P}}_{\alpha,a}^{f}(-\IF,\IF) }$ instead of $\ehe{{\mathcal{P}}_{\alpha,a}^{f}[0,\IF) }$. The same comments hold for all our results below.}
\end{remarks}
\COM{
\BK\label{Thm5}
Suppose that $\vk{X}(t)$ satisfies \eqref{cc2} and \eqref{lrr2}
and  $g(\vk{x})$ with order $p\in(0,\IF)$ satisfies \eqref{gne}.  Further, set $\alpha^*=\alpha p$,and $\beta^*=\beta p$ then we have
\BQNY
\pk{\sup_{t\in[0,T]}g(\vk{X}(t))>u}\sim u^{(\frac{2}{\alpha^*}-\frac{2}{\beta^*})_{+}}\pk{g(\vk{X}(t_0))>u}
\LT\{
\begin{array}{ll}
\mathbf{Q}a^{\frac{1}{\alpha}}\mathcal{H}_{\alpha}g^{-\frac{2}{\alpha p}}\int_{0}^\IF e^{-f(t)}dt,& \text{if}\ \beta^*>\alpha^*,\\
\widetilde{\mathcal{P}}_{\alpha,\frac{a}{g^{2/p}}}^{f}, & \text{if}\ \beta^*=\alpha^*,\\
1,& \text{if}\ \beta^*<\alpha^*,
\end{array}
\RT.
\EQNY
where $f(t)=\frac{b}{ g^{2/p}}t^\beta$.
\EK
}
\COM{
\BT\label{Thm4}
For a Gaussian vector process $\vk{X}(t)$, a continuous function $h(t)$ and a homogeneous function $g(\vk{x})$ with order $p\in[2,\IF)$, we assume that \eqref{gne}, \eqref{cc2}, and \eqref{lrr2} are satisfied.  Further, set $\alpha^*=\alpha p,\ \beta^*=\beta p$ and $f(t)=\frac{b\abs{t}^\beta}{g^{2/p}}$.\\
 \underline{(1)} If $p=2$, then we have
\BQNY
&&\pk{\sup_{t\in[0,T]}(g(\vk{X}(t))+h(t))>u}\\
&&\quad\quad\sim u^{(\frac{2}{\alpha^*}-\frac{2}{\beta^*})_{+}}\pk{g(\vk{X}(t_0))>u}
e^{\frac{h(t_0)}{2g}}
\LT\{
\begin{array}{ll}
\mathbf{Q}a^{\frac{1}{\alpha}}\mathcal{H}_{\alpha}g^{-\frac{2}{\alpha p}}\int_{0}^\IF e^{-f(t)}dt,& \text{if}\ \beta^*>\alpha^*,\\
\widetilde{\mathcal{P}}_{\alpha,\frac{a}{g^{2/p}}}^{f}, & \text{if}\ \beta^*=\alpha^*,\\
1,& \text{if}\ \beta^*<\alpha^*.
\end{array}
\RT.
\EQNY
\underline{(2)} If $p>2$, then we have
\BQNY
&&\pk{\sup_{t\in[0,T]}(g(\vk{X}(t))+h(t))>u}\\
&&\quad\quad\sim u^{(\frac{2}{\alpha^*}-\frac{2}{\beta^*})_{+}}\pk{g(\vk{X}(t_0))>u}
\LT\{
\begin{array}{ll}
\mathbf{Q}a^{\frac{1}{\alpha}}\mathcal{H}_{\alpha}g^{-\frac{2}{\alpha p}}\int_{0}^\IF e^{-f(t)}dt,& \text{if}\ \beta^*>\alpha^*,\\
\widetilde{\mathcal{P}}_{\alpha,\frac{a}{g^{2/p}}}^{f}, & \text{if}\ \beta^*=\alpha^*,\\
1,& \text{if}\ \beta^*<\alpha^*.
\end{array}
\RT.
\EQNY
\ET}
\subsection{Locally-stationary cases}
In this section, we assume that $X_i(t), i=1,\ldots, d$, are i.i.d. centered stationary Gaussian processes with unit variance function and covariance function $r(t)$ satisfying satisfying
\BQN\label{lrr11}
r(t,t+s)=1-a(t)|s|^{\alpha} +o(|s|^{\alpha}),\ \ s\rw 0,\ \alpha\in(0,2],
\EQN
where $a(t)$ are positive continuous function on $[0,T]$. Further assume
\begin{align}\label{lrr21}
r(t,t+s)<1,\ \  \forall\ t,t+s\in[0,T] \quad \hbox{and}\quad s\neq 0.
\end{align}
Write below $h_m$ for the maximum of a continuous function $h(t), t\in [0,T]$ and define
$$H:=\LT\{s\in[0,T]:h(s)=h_m\RT\}.$$	
	If $h_m$ is attained as some point $t_0 \in [0,T]$ it is important to know the behaviour of $h(t_0) -h(t)$ for $t$ close to $t_0$. Specifically, we shall assume that
\eqref{ht2} is satisfied for some $c> 0$, and $\gamma>0$.
\COM{	\BQN \label{eq:htt0}
	h_m - h(t) \sim c|t-t_0|^\gamma,\ \ \ \ \ t \rw t_0
	\EQN}

\COM{\BT\label{Thm0}
 Let $\vk{X}(t),t\in [0,T]$ be a Gaussian vector process with continuous sample paths and let $g(\vk{x})$ be a homogeneous function of order $p>0$. If \eqref{gne}, \eqref{lrr11} and \eqref{lrr21} are satisfied, then we have
\BQNY
\pk{\sup_{t\in[0,T]}Y(t)>u}
\sim \int_{0}^T(a(t))^{\frac{1}{\alpha}}dt\mathcal{H}_\alpha
u^{\frac{2}{\alpha p}}\pk{Y(0)>u},\quad u\rw \IF.
\EQNY
\ET
The form of $p_{h}(u), u\rw\IF$ depends on $p$, leading to two scenarios: $p\in(0,2)$ and $p\geq 2$. Next we give the results with these two scenarios.}

\BT\label{Thm1}
Suppose that Gaussian vector process $\vk{X}(t)$ with continuous sample paths satisfies \eqref{lrr11} and \eqref{lrr21}, a homogeneous function $g(\vk{x})$ with order $p>0$ satisfies \eqref{gne}, and $h(t)$ is a continuous function. Further, set $\alpha^*=\alpha p,$ and $ \beta^*=\frac{2\gamma p}{2-p}$.\\
\underline{(1)} If $h(t)\equiv0$, then  as $u\rw \IF$
\BQNY
\phu
\sim \mathcal{H}_\alpha \int_{0}^T(a(t))^{\frac{1}{\alpha}}dt
u^{\frac{2}{\alpha^*}}\pk{Y(0)>u}.
\EQNY

 \underline{(2)} For $p\in(0,2)$, if $H=\{t_0\}$ and \eqref{ht2} holds for some $c> 0, \gamma>0$, then we have as $u\rw \IF$
\BQNY
\phu \sim \mathbf{C}_{t_0}u^{(\frac{2}{\alpha^*}-\frac{2}{\beta^*})_{+}}\pk{Y(0)>u-h_m},
\EQNY
where $\mathbf{C}_{t_0}$ is the same as in \eqref{Ct0}
with $f(t)=\frac{c}{p}t^\gamma$ and $a=a(t_0)$.\\
\COM{\underline{(3)} If $H=\{t_1,\ldots, t_n\}\subset(0,T)$ with $t_1<t_2<\ldots<t_n$, $n\geq1$, and
\BQNY
h(t)=h_m-c_{j}|t-t_j|^{\gamma_{j}}(1+o(1)),\  t\rw t_j, j=1,\ldots n,
\EQNY
for some constants $\gamma_{j}>0$, $c_{j}>0$, then
\BQNY
&&\pk{\sup_{t\in[0,T]}(Y(t)+h(t))>u}\sim \eE{u^{(\frac{2}{\alpha^*}-\frac{2}{\beta^*})_+}} \pk{Y(0)>u-h_m} \\
&&\times \LT\{
\begin{array}{ll} \sum_{j=1}^n \mathbb{I}_{\{\gamma_j=\gamma\}}2
a_j^{1/\alpha}\mathcal{H}_{\alpha}g^{-\frac{2}{\alpha p}}\int_{0}^\IF e^{-f(t)}dt,&\ \ \text{if}\ \beta^*>\alpha^*,\\
\sum_{j=1}^n\LT(\mathbb{I}_{\{\gamma_j=\frac{\alpha(2-p)}{2}\}}
\mathcal{P}_{\alpha,
\frac{a_j}{g^{2/p}}}^{f}(-\IF,\IF)+
\mathbb{I}_{\{\gamma_j<\frac{\alpha(2-p)}{2}\}}\RT), &\ \ \text{if}\ \beta^*\leq\alpha^*,
\end{array}
\RT.
\EQNY
where $\gamma=\max_{1\leq j\leq n}\gamma_j$ and $a_j=a(t_j)$.\\}
\underline{(3)} For $p\in[2,\IF)$, then we have as $u\rw \IF$
\BQNY
\phu \sim  \mathcal{H}_\alpha \int_{0}^T (a(t))^{\frac{1}{\alpha}}e^{\frac{h(t) }{2}\mathbb{I}_{\{p=2\}}}dtu^{\frac{2}{\alpha^*}}\pk{Y(0)>u}.
\EQNY

\ET
\begin{remarks}
i) If $H$ consists of $n$ discrete points, say $t_1 \ldot t_n$, then as mentioned in \cite{Pit96} the tail of the supremum is easily obtained assuming that for each $t_i$ the assumptions of \netheo{Thm1} statement {i)} hold, implying that
$$ \phu 
\sim \Bigl(\sum_{j=1}^n \mathbf{C}_{t_j} \Bigr)  u^{(\frac{2}{\alpha^*}-\frac{2}{\beta^*})_+}\pk{Y(0)>u-h_m} .
$$
ii) If $H=[A,B]\subset[0,T]$ with $0\leq A<B\leq T$, then as $u\to \IF$
\BQN\label{intervalmax}
\phu 
\sim \mathcal{H}_\alpha  \int_{A}^B(a(t))^{1/\alpha}dtu^{\frac{2}{\alpha^*}}\pk{Y(0)>u-h_m},
\EQN
which is proven in Appendix.\\
iii) In \netheo{Thm1}, we investigate that if $p>2$, then $h$ does not contributed to the asymptotics.
\end{remarks}
We present next some examples of Gaussian chaos processes.
\cLb{\begin{exxa} \label{ex1}
 i) ($L^\rho$ norm process) For $g(\vk{x})=||\vk{x}||^p_\rho,\ p>0$ with
 and $$||\vk{x}||_\rho=\LT\{
\begin{array}{ll}
\LT(\sum_{i=1}^d|x_{i}|^\rho\RT)^{1/\rho},&\ \rho\in[1,\IF),\\
\max(|x_{1}|,\ldots,|x_{d}|),&\ \rho=\IF,
\end{array}
\RT.$$
we have \netheo{Thm3} and \netheo{Thm1} hold with order $p>0$.\\
 ii) (Product of several i.i.d Gaussian processes)
For $g(\vk{x})=\Pi_{i=1}^dx_i$, we have \netheo{Thm3} and \netheo{Thm1} hold with order $p=d$.\\
 iii) (Maximum of several i.i.d Gaussian processes)
For $g(\vk{x})=\max_{1\leq i\leq d}x_i$, we have \netheo{Thm3} and \netheo{Thm1} hold with order $p=1$.
\end{exxa}}

\section{Proofs}
We first give several preliminary lemmas, which play an important role in the proof of \netheo{Thm3} and \netheo{Thm1}. The proof of lemmas are shown in Appendix.
\ehe{ In the following we set $$\vk{ \xi}_{\vk{f}, u}(t)=(\xi_{1, f_1, u}(t) \ldot \xi_{d, f_d, u}(t))   $$
where $\xi_{i, f_i, u}(t) = \frac{ \xi_i( u^{-2/\alpha} t+ t_0)}{1+ u^{-2} f_i(t)} , $
where $f_i$'s given function and $\vk{\xi}(t)= (\xi_1(t) \ldot \xi_d(t)$ is a centeref Gaussian vector process.
}

\BEL\label{tt}
\ehe{Let} $\vk{X}(t)=(X_1(t)\ldots,X_d(t)),\ t\in[0,T],d\ge 1$ be a centered continuous  vector process with independent marginals which have unit variances and correlation functions satisfying \eqref{lrr21}. If $g: \R^d \to \R$ is a measurable function such that for some $p>0$
$$\abs{g(\vk x)} \le  \norm{\vk x}^p, \quad \forall \vk{x} \in \R^d,$$
 where $\norm{\cdot}$ is the Euclidean norm and then for $ 0<t_1<t_2<t_3$ and $u$ large enough
\BQNY
\pk{\sup_{t\in[0,t_1]} Y(t)>u,\sup_{t\in[t_2,t_3]} Y(t)>u}\leq
2\Psi\LT(\frac{2u^{\frac{1}{p}}-D}{\sqrt{4-\delta}}\RT),
\EQNY
where $D$ is a constant. Further, by \eqref{tdf} we have
\BQN\label{tt1}
\pk{\sup_{t\in[0,t_1]}Y(t)>u,\sup_{t\in[t_2,t_3]}
Y(t)>u}=o\LT(\pk{Y(0)>u}\RT), \ u\rw \IF.
\EQN
\EEL

\COM{
The next lemma is crucial for the asymptotics of supremum on short intervals. We shall assume that \eqref{aufrichten} holds with
$$ \lambda= 2/ p,\ L= \frac{1}{2g^{2/p}}. $$
This implies that $X(0)$ has distribution function in the Gumbel max-domain of attraction with scaling function $w(u)= L \lambda u^{\lambda-1}$,
i.e.,
\BQNY
\pk{ X(0)> u + s/w(u)} \sim e^{-s} \pk{X(0)> u}, \quad \forall s\inr.
\EQNY
We shall assume below a slightly stronger conditions on the pdf $f$ of  $X(0)$, namely
\BQN\label{strongF}
f(u) \sim  w(u) \pk{ X(0)> u} , \quad u\to \IF.
\EQN}

\BEL\label{im1}
let $\vk{\xi}(t)=(\xi_1(t)\ldots,\xi_d(t)),\ t\in\R,d\ge 1$ be a centered continuous  vector process with independent marginals which have unit variances and correlation functions satisfying \eqref{lrr11}.
Further we assume that $g(\vk{x}),\vk{x}\in\R^d, d\geq 1$, is a non-negative homogeneous function satisfying \necon{gcon} with order $1$.
Set $a=a(t_0), t_0\in\R$ and $K_u$ a family of index sets and $u_k$ satisfying that
\BQN\label{uk}
\lim_{u\rw\IF}\sup_{k\in K_u}\LT|\frac{u_k}{u}-1\RT|=0.
\EQN
If \ehe{$f_i,i\le d$ are  a continuous functions vanishing at 0}  , 
then we have that for some constants $S_1, S_2\geq 0$ and $S_1+ S_2>0$
\BQN\label{pp1}
\lim_{u\rw\IF}\sup_{k\in K_u}\LT|\frac{\pk{\sup_{t\in [-S_1,S_2]} \ehe{  g( \vk{ \xi}_{\vk{f}, u}(t)) }
			>u_k}}
{\pk{g(\vk{\xi}(t_0))>u_k}}-
\mathcal{P}_{\alpha,a}^{\vk f}
[-S_1,S_2]\RT|=0.
\EQN


\COM{and
\begin{align*}\label{pp2}
\mathcal{P}_{\alpha,a}^{f(t)}
[-S_1,S_2]=\E{\exp\LT(\sup_{t\in[-S_1,S_2]}
\sqrt{2a}B_{\alpha}(t)-a
\abs{t}^\alpha-f(t)\RT)}
.
\end{align*}}
If $\lim_{u\rw\IF}\sup_{k\in K_u}\abs{k u^{-2/\alpha}}\leq T$ 
some small enough $T>0$,
we have for some positive constant $S$ that when $u$ large enough
\BQN\label{P2p2}
\mathcal{H}_\alpha[0,(a-\vn_{T})^{1/\alpha}S]
&\leq&\frac{\pk{\sup_{t\in [0,S]}g(\vk{\xi}(u^{-2/\alpha}t+ku^{-2/\alpha}S+t_0))>u_k}}
{\pk{g(\vk{\xi}(t_0))>u_k}}\nonumber\\
&\leq& \mathcal{H}_\alpha[0,(a+\vn_{T})^{1/\alpha}S],
\EQN
holds for any $k\in K_u$ where $\vn_T\rw 0$, as $T\rw 0$.
Specially, if $T=0$, we have
\BQN\label{P2p}
\lim_{u\rw\IF}\sup_{k\in K_u}\LT|\frac{\pk{\sup_{t\in [0,S]}g(\vk{\xi}(u^{-2/\alpha}t+ku^{-2/\alpha}S+t_0))>u_k}}
{\pk{g(\vk{\xi}(t_0))>u_k}}-\mathcal{H}_\alpha[0,a^{1/\alpha}S]\RT|=0.
\EQN

\EEL

\BEL\label{in1}
Let the Gaussian vector process  $\vk{\xi}(t),t\in \R$ with independent marginals which have common correlation function $r(t)$ \cLb{satisfying \eqref{lrr11}} and the homogeneous function $g(\vk{x})$ satisfies \necon{gcon} with order $1$. Further, \cLb{ for some $t_0\in[0,T]$, set $a=a(t_0),$} and let  $K_u$ be  a family of \ehe{ countable} index sets such that for given positive constants  $u_k, k \in K_u$ we have
\BQN
\lim_{u\rw\IF}\sup_{k\in K_u}\LT|\frac{u_k}{u}-1\RT|=0.
\EQN
If  $\vn_0$ is  such that for all $t\in[t_0-\vn_0,t_0+\vn_0]$
\BQNY
\frac{a}{2}\abs{t-t_0}^{\alpha}\leq1-r(t)\leq2a\abs{t-t_0}^\alpha,
\EQNY
then we can find a constant $\mathbb{C}$ such that for all $S>0$ and $ T_2-T_1>S$ we have
\BQNY
{\limsup}_{u\rw\IF}\sup_{k\in K_u}\frac{\pk{\mathcal{A}_1(u_k),\mathcal{A}_2(u_k)}}
{\pk{g(\vk{\xi}(t_0))>u_k}}\leq \mathbb{C}\exp\LT(-\frac{\underline{a}}{8}|T_2-T_1-S|^\alpha\RT),
\EQNY
where $\mathcal{A}_i(u_k)=\LT\{\sup_{t\in [T_i,T_i+S]}g(\vk{\xi}(u^{-2/\alpha}(t+kS)+t_0))>u_k\RT\},\ i=1,2$, $\underline{a}=\inf_{t\in[t_0-\vn_0,t_0+\vn_0]}a(t)$  and
\BQN
\lim_{u\rw\IF}\sup_{k\in K_u}\abs{u^{-2/\alpha}k}\leq \vn_0.
\EQN
\EEL
\COM{In the next remark, since we couldn't make sure that $\sup_{k\in K_u}(u^{-2/\alpha^*}kS)\rw 0, u\rw \IF $, a more careful consideration is needed.
\begin{remark}
By \nelem{in1}, if we set $T_1=kS, T_2=(l+k)S$ with $k\in K_u, l>1,l\in\N$, then
\BQNY
&&{\limsup}_{u\rw\IF}\sup_{k\in K_u}\frac{\pk{\sup_{t\in [kS,(k+1)S]}g(\vk{\xi}(u^{-2/\alpha^*}t+t_0))>u_k,\sup_{t\in [(l+k)S,(l+k+1)S]}g(\vk{\xi}(u^{-2/\alpha^*}t+t_0))>u_k}}
{\pk{g(\vk{\xi}(t_0))>u_k}}\\
&&\leq{\limsup}_{u\rw\IF}\sup_{k'\in K_u}\sup_{k\in K_u}\frac{\pk{\sup_{t\in [k',(k'+1)S]}g(\vk{\xi}(u^{-2/\alpha^*}t+t_0))>u_k,\sup_{t\in [(l+k')S,(l+k'+1)S]}g(\vk{\xi}(u^{-2/\alpha^*}t+t_0))>u_k}}
{\pk{g(\vk{\xi}(t_0))>u_k}}\\
&&\leq \mathbb{C}\exp\LT(-\frac{a}{8}|(l-1)S|^\alpha\RT).
\EQNY
\end{remark}}

In the following proofs, $\mathbb{Q}_i,i\in\N$ denote some positive constants which can be different from line by line.
Further, in the proofs of \netheo{Thm3} and \netheo{Thm1}, we denote for some sets $\Delta_1, \Delta_2\subseteq\R$
\begin{align}
&\pu{\Delta_1}=\pk{\sup_{\Delta_1}(Y(t)+h(t))>u},\quad \pu{\Delta_1,\Delta_2}=\pk{\sup_{\Delta_1}(Y(t)+h(t))>u,\sup_{\Delta_2}(Y(t)+h(t))>u},
\label{ppu}\\
&\qu{\Delta_1}=\pk{\sup_{\Delta_1}Y(t)>u},\quad \qu{\Delta_1,\Delta_2}=\pk{\sup_{\Delta_1}Y(t)>u,\sup_{\Delta_2}Y(t)>u}.\nonumber
\end{align}
\prooftheo{Thm3}
We present first the proof for  $t_0=0$. Without loss of generalities we shall assume that $h(0)=0$.\\
For all $u$ large, let $\Delta(u)=[0,\delta(u)]$, where
$$\delta(u)=\LT(\frac{(\ln u)^q }{u}\RT)^{2/\beta^*}, \quad q>\max(\frac{p}{2},\frac{p}{2-p}).$$
 It follows that for $\theta>0$ small enough
\BQN\label{Thm211}
\pu{\Delta(u)}\le \pu{[0,T]}\le \pu{\Delta(u)}+ \pu{[\delta_u,\theta]}+ \pu{[\theta,T]}.
\EQN
We  first give upper bounds of $\pu{[\delta_u,\theta]}$ and $\pu{[\theta,T]}$ which finally imply that as $u\to \IF$
\BQN\label{Thmeq02}
\pu{[\delta_u,\theta]}=o(\pu{\Delta(u)}),\ \ \pu{[\theta,T]}=o(\pu{\Delta(u)}) .
\EQN
Set $\overline{\vk{X}}(t)=(\overline{X}_1(t)  \ldot \overline{X}_d(t))$ with
$\overline{X}_i(t)=\frac{X_i(t)}{\sigma(t)}, i\le d$ and define
$$h^*=\max_{t\in[0,T]}h(t), \  \sigma_\theta=\sup_{t\in[\theta,T]}\sigma(t)<1.$$ Then by Borell inequality
\BQN\label{PI2}
\pu{[\theta,T]}&\leq& \pk{\sup_{t\in[\theta,T]} Y(t)>u-h^*}\leq \pk{\sup_{t\in[\theta,T]} g(\overline{\vk{X}}(t))>\frac{u-h^*}{\sigma_\theta^p}}\nonumber\\
&\leq& \pk{\sup_{t\in[\theta,T]}\sup_{\vk{z}\in\mathbb{S}_{d-1}} \langle\overline{\vk{X}}(t),\vk{z}\rangle>\frac{(u-h^*)^{1/p}}
{\sigma_\theta}}\leq \exp\LT(-\frac{\LT(\frac{(u-h^*)^{1/p}}{\sigma_\theta}
-\mathbb{Q}_0\RT)^2}{2}\RT)\nonumber\\
&=&o\LT(\pk{Y(0)>u}\RT),\ u\rw\IF,
\EQN
where $\mathbb{Q}_0:=\E{\sup_{t\in[\theta,T],\vk{z}\in\mathbb{S}_{d-1}}
\langle\overline{\vk{X}}(t),\vk{z}\rangle}$.
By \eqref{cc2} and \eqref{ht2}, we know that for some $\vn,\vn_1 \in(0,1)$
\BQN
&\frac{(u-h(t))^{1/p}}{\sigma(t)}
\geq u^{1/p}(1+\frac{c(1-\vn)}{up}|t|^\gamma)(1+(1-\vn)b|t|^\beta)
\geq u^{1/p}\LT(1+\frac{c(1-\vn)}{up}|t|^\gamma+(1-\vn) b|t|^\beta\RT),\label{hl2}\\
&\frac{(u-h(t))^{1/p}}{\sigma(t)}
\leq u^{1/p}(1+\frac{c(1+\vn_1)}{up}|t|^\gamma)(1+(1+\vn_1)b|t|^\beta)
\leq u^{1/p}\LT(1+\frac{c(1+\vn)}{up}|t|^\gamma+(1+\vn) b|t|^\beta\RT)\label{hu2}
\EQN
holds for $t\in[0,\theta]$, then
\begin{align*}
\inf_{t\in[\delta(u),\theta] }\frac{(u-h(t))^{2/p}}{\sigma^2(t)}
&\geq \inf_{t\in[\delta(u),\theta] } u^{2/p}\LT(1+\frac{c(1-\vn)}{up}|t|^\gamma+(1-\vn) b|t|^\beta\RT)^{2}\\
&\geq u^{2/p}+\mathbb{Q}_1(\ln u)^{\frac{q(2-p)}{p}\vee\frac{2q}{p}}.
\end{align*}
By \eqref{lrr2}, we have that
\begin{align*}
\E{\LT(\langle\overline{\vk{X}}(t),\vk{z}\rangle\RT)^2}=1
\end{align*}
and
\begin{align*}
\E{\LT(\langle\overline{\vk{X}}(t),\vk{z}\rangle-\langle
\overline{\vk{X}}(s),\vk{z}'\rangle\RT)^2}
&\leq 2\E{\LT(\langle\overline{\vk{X}}(t),\vk{z}\rangle-\langle
\overline{\vk{X}}(s),\vk{z}\rangle\RT)^2}
+2\E{\LT(\langle\overline{\vk{X}}(s),\vk{z}\rangle-\langle
\overline{\vk{X}}(s),\vk{z}'\rangle\RT)^2}\\
&\leq 2\E{\LT(\langle\overline{\vk{X}}(t)-\overline{\vk{X}}(s),\vk{z}\rangle\RT)^2}
+2\E{\LT(\langle\overline{\vk{X}}(s),\vk{z}-\vk{z}'\rangle\RT)^2}\\
&\leq \mathbb{Q}_2|s-t|^\alpha+ \mathbb{Q}_3\sum_{i=1}^d|z_i-z'_i|^2
\leq \mathbb{Q}_4\LT(|s-t|^\alpha+\sum_{i=1}^d|z_i-z'_i|^\alpha\RT)
\end{align*}
holds for $s,t\in[0,\theta]$ and $\vk{z},\vk{z}'\in\mathbb{S}_{d-1}$.
Then it follows from \cite{Pit96}[Theorem 8.1] that
\begin{align*}
\pu{[\delta_u,\theta]}&\leq \pk{\sup_{t\in[\delta(u),\theta]} g(\overline{\vk{X}}(t))>\inf_{t\in[\delta(u),\theta]}\frac{u-h(t)}{\sigma^p(t)}}\leq \pk{\sup_{t\in[\delta(u),\theta]}\sup_{\vk{z}\in\mathbb{S}_{d-1}} \langle\overline{\vk{X}}(t),\vk{z}\rangle>\inf_{t\in[\delta(u),\theta]}
\frac{(u-h(t))^{1/p}}{\sigma(t)}}\\
&\leq\mathbb{Q}_5u^{\frac{2(d+1)}{\alpha}}\Psi\LT(\inf_{t\in[\delta(u),\theta]}
\frac{(u-h(t))^{1/p}}{\sigma(t)}\RT)
=o\LT(\pk{Y(0)>u}\RT),\ u\rw\IF,
\end{align*}
where  in the last equation we use the fact that
\BQNY
\inf_{t\in[\delta(u),\theta] }\frac{(u-h(t))^{2/p}}{\sigma^2(t)}-u^{2/p}
\geq \mathbb{Q}_1(\ln u)^{\frac{q(2-p)}{p}\vee\frac{2q}{p}}
\rw \IF, \ u\rw\IF.
\EQNY
Next, we give the asymptotic of $\pu{\Delta(u)}$, as $u\rw\IF$.
Set for any $S>0$,
\BQNY
&&I_k(u)=[ku^{-2/\alpha^*}S,(k+1)u^{-2/\alpha^*}S],\ \  k\in\N,\ \  N(u)=\LT\lfloor(\ln u)^{\frac{2q}{\beta^*}}u^{\frac{2}{\alpha^*}-\frac{2}{\beta^*}}S^{-1}\RT\rfloor.\\
&&\mathcal{G}_{u,+\vn}(k)=u^{1/p}\LT(1+\frac{c(1+\vn)}{up}
\abs{(k+1)u^{-2/\alpha^*}S}^\gamma+(1+\vn) b\abs{(k+1)u^{-2/\alpha^*}S}^\beta\RT), \\
&&\mathcal{G}_{u,-\vn}(k)=u^{1/p}\LT(1+\frac{c(1-\vn)}{up}
\abs{ku^{-2/\alpha^*}S}^\gamma+(1-\vn) b\abs{ku^{-2/\alpha^*}S}^\beta\RT).
\EQNY
{\bf Case 1}: $\beta^*>\alpha^*$. For $u$ large enough we have
\BQN
\sum_{k=0}^{N(u)}\pu{I_{k}(u)}\geq \pu{\Delta(u)}\geq\sum_{k=0}^{N(u)-1}\pu{I_{k}(u)}-\sum_{i=1}^2\Lambda_i(u),\label{Thmeq12}
\EQN
where
$$\Lambda_1(u)=\sum_{k=0}^{N(u)}\pu{ I_{k}(u),I_{k+1}(u)},\quad \Lambda_2(u)=\sum_{0\leq k,l\leq N(u), l\geq k+2}\pu{I_{k}(u),I_{l}(u)}.$$
 Set $\widetilde{g}(\vk{x})=g^{1/p}(\vk{x})$, then $\widetilde{g}(\vk{x})$ is homogeneous function with order $1$ and satisfies \necon{gcon}.\\
Then in light of \nelem{im1} and \eqref{hl2}, we have that for some $\epsilon\in [0,1)$,
\begin{align}\label{Th21}
\sum_{k=0}^{N(u)}\pu{I_{k}(u)}\leq& \sum_{k=0}^{N(u)}\pk{\sup_{t\in I_k(u)}\widetilde{g}(\overline{\vk{X}}(t))>\mathcal{G}_{u,-\vn}(k)}\nonumber\\
=& \sum_{k=0}^{N(u)}\pk{\sup_{t\in [0,S]}\widetilde{g}(\overline{\vk{X}}(u^{-2/\alpha^*}t+ku^{-2/\alpha^*}S))>\mathcal{G}_{u,-\vn}(k)}\nonumber\\
\sim& h_0\mathcal{H}_\alpha[0,a^{1/\alpha}S]u^{\frac{m-1}{p}}
\sum_{k=0}^{N(u)}\exp\LT(-\frac{1}{2}\LT(\mathcal{G}_{u,-\vn}(k)\RT)
^{2}\RT)\nonumber\\
\sim& h_0\mathcal{H}_\alpha[0,a^{1/\alpha}S]u^{\frac{m-1}{p}}
\exp\LT(-\frac{u^{\frac{2}{p}}}{2}\RT)\nonumber\\
&\times
\sum_{k=0}^{N(u)}\exp\LT(-\frac{c}{p}(1-\vn-\epsilon)u^{\frac{2-p}{p}}|k S u^{-\frac{2}{\alpha^*}}|^\gamma-b(1-\vn-\epsilon)u^{2/p}|k S u^{-\frac{2}{\alpha^*}}|^\beta\RT)\nonumber\\
\sim& h_0\mathcal{H}_\alpha[0,a^{1/\alpha}S]u^{\frac{m-1}{p}}
\exp\LT(-\frac{u^{\frac{2}{p}}}{2}\RT)
\sum_{k=0}^{N(u)}\exp\LT(-(1-\vn-\epsilon)f(u^{\frac{2}{\beta^*}}k S u^{-\frac{2}{\alpha^*}})\RT)\nonumber\\
\sim&
\pk{Y(0)>u}\frac{\mathcal{H}_\alpha[0,a^{1/\alpha}S]}{S}
u^{\frac{2}{\alpha^*}-\frac{2}{\beta^*}}
\int_0^{\IF}\exp\LT(-(1-\vn-\epsilon)f(t)\RT)dt\nonumber\\
\sim&
\pk{Y(0)>u}a^{1/\alpha}\mathcal{H}_\alpha u^{\frac{2}{\alpha^*}-\frac{2}{\beta^*}}
\int_0^{\IF}e^{-f(t)}dx,
\end{align}
as $ u\rw\IF, \ S\rw\IF,\ \vn\rw 0, \ \epsilon\rw 0$ where
$f(t)=\frac{c}{p}\abs{t}^\gamma\mathbb{I}_{\{\beta^*=\frac{2\gamma p}{2-p}\}} +b\abs{t}^\beta\mathbb{I}_{\{\beta^*=\beta p\}} $.\\
 Similarly, we derive that as $u\rw\IF, \ S\rw\IF$,
\BQN\label{Th212}
\sum_{k=0}^{N(u)-1}\pu{I_{k}(u)}\geq\pk{Y(0)>u}a^{1/\alpha}\mathcal{H}_\alpha u^{\frac{2}{\alpha^*}-\frac{2}{\beta^*}}
\int_0^{\IF}e^{-f(t)}dt.
\EQN
Moreover,
\begin{align}\label{Thmeq22}
\Lambda_1(u)\leq&\sum_{k=0}^{N(u)}\left(\pu{I_{k}(u)}+\pu{I_{k+1}(u)}-\pu{(I_{k}(u)\cup I_{k+1}(u))}\right)\nonumber\\
\leq& \sum_{k=0}^{N(u)}\left(\pk{ \sup_{t\in I_{k}(u)} \widetilde{g}(\overline{\vk{X}}(t))>\mathcal{G}_{u,-\vn}(k)}+\pk{\sup_{t\in  I_{k+1}(u)} \widetilde{g}(\overline{\vk{X}}(t))>\mathcal{G}_{u,-\vn}(k)}\RT.\nonumber\\
 &\LT.-\pk{\sup_{t\in ((I_{k}(u)\cup I_{k+1}(u)))} \widetilde{g}(\overline{\vk{X}}(t))>\widehat{\mathcal{G}}_{u,-\vn}(k)}\right)\nonumber\\
\leq&   h_0\LT(2\mathcal{H}_\alpha[0,a^{1/\alpha}S]-
\mathcal{H}_\alpha[0,2a^{1/\alpha}S]\RT)
u^{\frac{m-1}{p}}
\sum_{k=0}^{N(u)}\exp\LT(-\frac{1}{2}
\LT(\widehat{\mathcal{G}}_{u,-\vn}(k)\RT)^{2}\RT)\nonumber\\
\sim&
\frac{2\mathcal{H}_\alpha[0,a^{1/\alpha}S]-
\mathcal{H}_\alpha[0,2a^{1/\alpha}S]}{S}
\int_0^{\IF}\exp\LT(-(1-\vn-\epsilon)f(t)\RT)dt u^{\frac{2}{\alpha^*}-\frac{2}{\beta^*}}\pk{Y(0)>u}\nonumber\\
=&o\left(u^{\frac{2}{\alpha^*}-\frac{2}{\beta^*}}\pk{Y(0)>u}\right), \ u\rw\IF, S\rw \IF, \vn\rw 0,\epsilon\rw 0.
\end{align}
where $\widehat{\mathcal{G}}_{u,-\vn}(k)=\min (\mathcal{G}_{u,-\vn}(k),\mathcal{G}_{u,-\vn}(k+1))$.
By \nelem{in1}, we have
\begin{align}\label{Thmeq32}
\Lambda_2(u)
&\leq\sum_{0\leq k,l\leq N(u), l\geq k+2}\pk{ \sup_{t\in I_{k}(u)}\widetilde{g}(\overline{\vk{X}}(t))>\mathcal{G}_{u,-\vn}(k), \sup_{t\in I_{l}(u)} \widetilde{g}(\overline{\vk{X}}(t))>\mathcal{G}_{u,-\vn}(l)}\nonumber\\
&\leq\sum_{0\leq k\leq N(u)}\sum_{l=2}^{N(u)}\pk{ \sup_{t\in I_{k}(u)}\widetilde{g}(\overline{\vk{X}}(t))>\mathcal{G}_{u,-\vn}(k), \sup_{t\in I_{k+l}(u)} \widetilde{g}(\overline{\vk{X}}(t))>\mathcal{G}_{u,-\vn}(k)}\nonumber\\
&\leq\mathbb{Q}_6\LT(\sum_{k=0}^{N(u)}\pk{\widetilde{g}(\overline{\vk{X}}(0))>
\mathcal{G}_{u,-\vn}(k)}
\RT)\sum_{l=2}^{\IF}\exp\LT(-(l S)^\alpha/8\RT)\nonumber\\
&\leq\mathbb{Q}_7\pk{Y(0)>u} u^{\frac{2}{\alpha^*}-\frac{2}{\beta^*}}S\sum_{l=0}^{\IF}\exp\LT(-(l S)^\alpha/8\RT)\nonumber\\
&= o\LT(u^{\frac{2}{\alpha^*}-\frac{2}{\beta^*}}\pk{Y(0)>u}\RT),\ u\rw\IF,\ S\rw\IF.
\end{align}
Combing (\ref{Th21})-(\ref{Thmeq32}) with (\ref{Thmeq12}), we obtain
\BQN\label{Th2re1}
\pu{\Delta(u)}\sim \pk{Y(0)>u}a^{1/\alpha}\mathcal{H}_\alpha u^{\frac{2}{\alpha^*}-\frac{2}{\beta^*}}
\int_0^{\IF}\exp\LT(-f(t)\RT)dt,\ u\rw\IF.
\EQN
{\bf Case 2:} $\beta^*=\alpha^*$.
We consider that for $u$ large enough
\BQN
\pu{I_0(u)}\leq \pu{\Delta(u)}\leq\sum_{k=0}^{N(u)}\pu{I_k(u)}.\label{Th22}
\EQN
Using Lemma \ref{im1} and \eqref{hu2}, we have that for some small $\epsilon\in(0,1)$
\begin{align}\label{Thmeq62}
\pu{I_0(u)}
&\geq \pk{\sup_{t\in [0,S]} \LT(Y(t u^{-2/\alpha^*})+h(t u^{-2/\alpha^*})\RT)>u}\nonumber\\
&\geq \pk{\sup_{t\in [0,S]} \frac{\widetilde{g}(\overline{\vk{X}}(t u^{-2/\alpha^*}))}{1+\frac{c(1+\vn)}{up}\abs{tu^{-2/\alpha^*}}^\gamma+(1+\vn) b\abs{tu^{-2/\alpha^*}}^\beta}>u^{1/p}}\nonumber\\
&\geq \pk{\sup_{t\in [0,S]} \frac{\widetilde{g}(\overline{\vk{X}}(t u^{-2/\alpha^*}))}{1+(1+\vn+\epsilon)u^{-2/p} f(t)}>u^{1/p}}\nonumber\\
&\sim\E{\sup_{t\in[0,S]}\exp\LT(\sqrt{2a}
B_{\alpha}(t)-a\abs{t}^{\alpha}-(1+\vn+\epsilon)f(t)\RT)}
\pk{Y(0)>u}\nonumber\\
&\sim\mathcal{P}_{\alpha,a}^{f}[0,\IF)
\pk{Y(0)>u},\ u\rw\IF, \vn\rw 0, \epsilon\rw 0, S\rw\IF.
\end{align}
Moreover, by \nelem{im1}
\begin{align}\label{Thmeq82}
{\sum_{k=1}^{N(u)}}\pu{I_{k}(u)}
&\leq {\sum_{k=1}^{N(u)}} \pk{\sup_{t\in I_k(u)}\widetilde{g}(\overline{\vk{X}}(t))>\mathcal{G}_{u,-\vn}(k)}\nonumber\\
&\sim h_0\mathcal{H}_\alpha[0,a^{1/\alpha}S]u^{\frac{m-1}{p}}
\sum_{k=1}^{N(u)}\exp\LT(-\frac{1}{2}
\LT(\mathcal{G}_{u,-\vn}(k)\RT)^{2}\RT)\nonumber\\
&\leq  h_0\mathcal{H}_\alpha[0,a^{1/\alpha}S]u^{\frac{m-1}{p}}
\exp\LT(-\frac{u^{\frac{2}{p}}}{2}\RT)\nonumber\\
&\quad\times
\sum_{k=1}^{N(u)}\exp\LT(-(1-\vn-\epsilon)f\LT(u^{\frac{2}{\beta^*}}k S u^{-\frac{2}{\alpha^*}}\RT)\RT)\nonumber\\
&\leq
\pk{Y(0)>u}\mathcal{H}_\alpha[0,a^{1/\alpha}S]
\sum_{k=1}^{\IF}\exp\LT(-\mathbb{Q}_8(k S)^\gamma\RT)\nonumber\\
&\sim
\mathbb{Q}_9\pk{Y(0)>u}\mathcal{H}_\alpha a^{1/\alpha}
S\exp\LT(-\mathbb{Q}_{10} S^\gamma\RT)\nonumber\\
&=o\LT(\pk{Y(0)>u}\RT),\ u\rw\IF,\ S\rw\IF.
\end{align}
Inserting (\ref{Thmeq82}) and  (\ref{Thmeq62}) into \eqref{Th22}, we have
\BQN\label{Th2re2}
\pu{\Delta(u)}\sim\mathcal{P}_{\alpha,a}^{f}[0,\IF)
\pk{Y(0)>u}, u\rw\IF.
\EQN
{\bf Case 3:} $\beta^*<\alpha^*$.
\BQN\label{Thmeq92}
\pu{\Delta(u)}\geq\pk{ Y(0)>u}.
\EQN
For any $\vn_2\in (0,1)$,  $\Delta(u)\subseteq [0,u^{-2/\alpha^*}\vn_2]$ when $u$ large enough. Then as $u\rw\IF,\ \vn_2\rw0$
\begin{align*}
\pu{\Delta(u)}\leq \pk{\sup_{t\in [0,u^{-2/\alpha^*}\vn_2]}\widetilde{g}(\overline{\vk{X}}(t))>u^{1/p}}
\sim\mathcal{H}_{\alpha}[0,a^{1/\alpha}\vn_2]\pk{Y(0)>u}\sim\pk{Y(0)>u}.
\end{align*}
Together with (\ref{Thmeq92}), we get
\BQN\label{Th2re3}
\pu{\Delta(u)}\sim \pk{ Y(0)>u}.
\EQN
Further, \eqref{Thmeq02} are derived according to \eqref{Th2re1}-\eqref{Th2re3}.\\
\ehe{Finally, we note that if  $t_0\in(0,T)$ and $t_0=T$, the proof is the same with as above by simply  replacing $\Delta(u)$ by  $\Delta(u)=[-\delta(u),\delta(u)]$ and $\Delta(u)=[-\delta(u),0]$, respectively.}
Thus we complete the proof.
\QED
\COM{
\proofkorr{Thm5} When $p\in(0,2)$, by \netheo{Thm3} the results are clear if we notice that
$\beta^*=\beta p$ and $h(t_0)=0$ and $c=0$.\\
When $p\in[2,\IF)$, in the proof of \netheo{Thm3}, if we take $\beta^*=\beta p,\ h(t)\equiv0,\ c=0$ and
$f(t)=\frac{b \abs{t}^\beta}{ \gg^{2/p}}$, then all argumentations  in the proof of \netheo{Thm3} still hold.
\QED}
\COM{
\prooftheo{Thm4}  \underline{(1)} We consider the case that $t_0=0$.
 We have
 \BQN\label{Thm211}
\Pi_1(u)\le \pk{\sup_{t\in[0,T]}(g(\vk{X}(t))+h(t))>u}\le \Pi_1(u)+ \Pi_2(u)
\EQN
 with
\begin{align*}
&\Pi_1(u):=\pk{\sup_{t\in[0,\theta]} (g(\vk{X}(t))+h(t))>u },\\
&\Pi_2(u):=\pk{\sup_{t\in[\theta,T]} (g(\vk{X}(t))+h(t))>u }.
\end{align*}
By \eqref{PI2}, we know $\Pi_2(u)=o\LT(\pk{g(\vk{X}(0))>u-h(0)}\RT),\ u\rw\IF.$
Set $M^1_\theta=\sup_{t\in[0,\theta]}h(t)$ and $M^2_\theta=\inf_{t\in[0,\theta]}h(t)$.
\eE{Then by \neremark{Thm5}}
\begin{align*}
\Pi_1(u)&\leq \pk{\sup_{t\in[0,\theta]} g(\vk{X}(t))>u-M^1_\theta }\\
        &\sim \mathbf{C}u^{(\frac{2}{\alpha^*}-\frac{2}{\beta^*})_{+}}
        \pk{g(\vk{X}(0))>u-M^1_\theta}\\
        &\sim\mathbf{C}u^{(\frac{2}{\alpha^*}-\frac{2}{\beta^*})_{+}}\pk{g(\vk{X}(0))>u}
        e^{\frac{M^1_\theta}{2g}}\\
        &\sim\mathbf{C}u^{(\frac{2}{\alpha^*}-\frac{2}{\beta^*})_{+}}\pk{g(\vk{X}(0))>u}
        e^{\frac{h(0)}{2g}},\ u\rw\IF,\ \theta\rw0.
\end{align*}
Similarly,
\begin{align*}
\Pi_1(u)\geq \pk{\sup_{t\in[0,\theta]} g(\vk{X}(t))>u-M^2_\theta }
        \sim\mathbf{C}u^{(\frac{2}{\alpha^*}-\frac{2}{\beta^*})_{+}}\pk{g(\vk{X}(0))>u}
        e^{\frac{h(0)}{2g}},\ u\rw\IF,\ \theta\rw0.
\end{align*}
Then then claims follow.\\
\underline{(2)} When $p>2$, if we notice that $$\pk{g(\vk{X}(0))>u+\mathbb{Q}}\sim\pk{g(\vk{X}(0))>u}$$.
Then by \neremark{Thm5} and
 \BQNY
&& \pk{\sup_{t\in[0,T]}(g(\vk{X}(t))+h(t))>u}\geq \pk{\sup_{t\in[0,T]}g(\vk{X}(t))>u-\sup_{s\in[0,T]}h(s)},\\
&& \pk{\sup_{t\in[0,T]}(g(\vk{X}(t))+h(t))>u}\leq \pk{\sup_{t\in[0,T]}g(\vk{X}(t))>u-\inf_{s\in[0,T]}h(s)},
\EQNY
the results are clear.}

\prooftheo{Thm1} \underline{(1)}
For any $\theta>0$ and $S>0$, set $\alpha^*=\alpha p$
\BQNY
&&I_k(\theta)=[k\theta,(k+1)\theta],\quad a_k=a(k\theta)\quad k\in \N,\quad N(\theta)=\LT\lfloor\frac{T}{\theta}\RT\rfloor,\\
&&J^k_l(u)=\LT[k\theta+lu^{-2/\alpha^*}S,k\theta+(l+1)u^{-2/\alpha^*}S\RT],\quad M(u)=\LT\lfloor\frac{\theta u^{2/\alpha^*}}{S}\RT\rfloor.
\EQNY
We have
\BQNY
\sum_{k=0}^{N(\theta)-1}\LT(\sum_{l=0}^{M(u)-1}\qu{J^k_l(u)}\RT)-\sum_{i=1}^{4}\mathcal{A}_i(u)\leq\qu{[0,T]}\leq\sum_{k=0}^{N(\theta)}
\qu{I_k(\theta)}\leq\sum_{k=0}^{N(\theta)}\LT( \sum_{l=0}^{M(u)}\qu{J^k_l(u)}\RT),
\EQNY
where
\BQNY
\mathcal{A}_i(u)=\sum_{(k_1,l_1,k_2,l_2)\in \mathcal{L}_i}\qu{J^{k_1}_{l_1}(u),J^{k_2}_{l_2}(u)},\ i=1,2,3,4,
\EQNY
with
\BQNY
&&\mathcal{L}_{1}=\LT\{0\leq k_1= k_2\leq N(\theta)-1,0\leq l_1+1=l_2\leq M(u)-1\RT\},\\
&&\mathcal{L}_{2}=\LT\{0\leq k_1+1= k_2\leq N(\theta)-1, l_1=M(u), l_2=0\RT\},\\
&&\mathcal{L}_3=\LT\{0\leq k_1+1<k_2\leq N(\theta)-1,0\leq l_1,l_2\leq M(u)-1\RT\},\\
&&\mathcal{L}_4=\LT\{0\leq k_1\leq k_2\leq N(\theta)-1,k_2-k_1\leq 1,0\leq l_1,l_2\leq M(u)-1\RT\}\setminus\LT(\mathcal{L}_1\cup \mathcal{L}_2\RT).
\EQNY
By \nelem{im1},
\begin{align*}
\sum_{k=0}^{N(\theta)}\LT( \sum_{l=0}^{M(u)}\qu{J^k_l(u)}\RT)
&=\sum_{k=0}^{N(\theta)}\LT(\sum_{l=0}^{M(u)}\pk{\sup_{t\in [0,S]}Y(k\theta+lu^{-2/\alpha^*}S+u^{-2/\alpha^*} t)>u}\RT) \\
&=\sum_{k=0}^{N(\theta)}\LT(\sum_{l=0}^{M(u)}\pk{\sup_{t\in [0,S]}\widetilde{g}\LT(\vk{X}(k\theta+lu^{-2/\alpha^*}S+u^{-2/\alpha^*} t)\RT)>u^{1/p}}\RT) \\
&\leq \sum_{k=0}^{N(\theta)}\LT(\sum_{l=0}^{M(u)}(a_k+\vn_\theta)^{\frac{1}{\alpha}}
\mathcal{H}_\alpha S \pk{Y(0)>u}\RT)\\
&\sim \LT(\sum_{k=0}^{N(\theta)}(a_k+\vn_\theta\theta)^{\frac{1}{\alpha}}\RT)
\mathcal{H}_\alpha u^{2/\alpha^*} \pk{Y(0)>u}\\
&\sim \int_{0}^T(a(t))^{1/\alpha}d tu^{-2/\alpha^*}\mathcal{H}_\alpha \pk{Y(0)>u}, \ u\rw\IF,\ S\rw\IF,\ \theta\rw 0.
\end{align*}
where $\widetilde{g}(\vk{x})=g^{1/p}(\vk{x})$.
Similarly,
\BQNY
\sum_{k=0}^{N(\theta)-1}\LT(\sum_{l=0}^{M(u)-1}\qu{J^k_l(u)}\RT)
\geq \int_{0}^T(a(t))^{1/\alpha}d tu^{-2/\alpha^*}\mathcal{H}_\alpha \pk{Y(0)>u}, \ u\rw\IF,\ S\rw\IF,\ \theta\rw 0.
\EQNY
Now we find an upper bound for $\sum_{i=1}^{4}\mathcal{A}_i(u)$,  by \nelem{im1}
\begin{align*}
\mathcal{A}_1(u)
&=\sum_{k=0}^{N(\theta)-1}\LT(\sum_{l=0}^{M(u)-1}\left(\qu{J^k_l(u)}+\qu{J^k_{l+1}(u)}
-\qu{J^k_l(u)\cup J^k_{l+1}(u)}\right)\RT)\nonumber\\
&\sim   \sum_{k=0}^{N(\theta)-1}\LT(\LT(\mathcal{H}_\alpha[0,(a_k+\vn_\theta)^
{\frac{1}{\alpha}}S]
+\mathcal{H}_\alpha[0,(a_k+\vn_\theta)^{\frac{1}{\alpha}}S]
-\mathcal{H}_\alpha[0,2(a_k-\vn_\theta)^{\frac{1}{\alpha}}S]\RT)
\sum_{l=0}^{M(u)-1}\pk{Y(0)>u}\RT)\nonumber\\
&\leq\mathbb{Q}_1\LT(\sum_{k=0}^{N(\theta)-1}\LT((a_k+\vn_\theta)^{\frac{1}{\alpha}}-(a_k
-\vn_\theta)^{\frac{1}{\alpha}}\RT)\theta\RT)u^{2/\alpha^*}\pk{Y(0)>u}\\
&=o\left(u^{2/\alpha^*}\pk{Y(0)>u}\right), \ u\rw\IF,\ \ S\rw\IF,\theta\rw0.
\end{align*}
Similarly,
\begin{align*}
\mathcal{A}_2(u)
&=\sum_{k=0}^{N(\theta)-1}\pu{J^{k}_{M(u)-1}(u),J^{k+1}_{0}(u)}\\
&\leq\sum_{k=0}^{N(\theta)-1}\pk{\sup_{t\in [0,2S] }Y((k+1)\theta-u^{-2/\alpha^*}t)>u,\sup_{t\in [0,2S]}Y((k+1)\theta+u^{-2/\alpha^*}t)>u}\\
&=\sum_{k=0}^{N(\theta)-1}\LT(\pk{\sup_{t\in [0,2S] }Y((k+1)\theta-u^{-2/\alpha^*}t)>u}+\pk{\sup_{t\in [0,2S]}Y((k+1)\theta+u^{-2/\alpha^*}t)>u}\RT.\\
&\quad\LT.-\pk{\sup_{t\in [-2S,2S] }Y((k+1)\theta-u^{-2/\alpha^*}t)>u}\RT)\nonumber\\
&\sim   \sum_{k=0}^{N(\theta)-1}\LT(\LT(2\mathcal{H}_\alpha[0,2(a_k+\vn_\theta)
^{\frac{1}{\alpha}}S]
-\mathcal{H}_\alpha[-2(a_k-\vn_\theta)^{\frac{1}{\alpha}}
S,2(a_k-\vn_\theta)^{\frac{1}{\alpha}}S]\RT)\sum_{l=0}^{M(u)-1}\pk{Y(0)>u}\RT)\nonumber\\
&\leq\mathbb{Q}_2\LT(\sum_{k=0}^{N(\theta)-1}\LT((a_k+\vn_\theta)^
{\frac{1}{\alpha}}-(a_k
-\vn_\theta)^{\frac{1}{\alpha}}\RT)\theta\RT)u^{2/\alpha^*}\pk{Y(0)>u}\\
&=o\left(u^{2/\alpha^*}\pk{Y(0)>u}\right), \ u\rw\IF,\ \ S\rw\IF,\theta\rw0.
\end{align*}

For any $\theta>0$
\BQNY
\E{X_i(t)X_i(s)}=r(s,t)\leq 1-\delta(\theta)
\EQNY
for $(s,t)\in J^{k_1}_{l_1}(u)\times J^{k_2}_{l_2}(u),(j_1,k_1,j_2,k_2)\in\mathcal{L}_3$ where $\delta(\theta)>0$ is related to $\theta$. Then by \nelem{tt} as $u\rw\IF, S\rw\IF, \theta\rw 0$
\begin{align*}
\mathcal{A}_3(u)\leq N(\theta)M(u)2\Psi\LT(\frac{2u^{\frac{1}{2}}
-\mathbb{Q}_3}{\sqrt{4-\delta(\theta)}}\RT)
\leq\frac{T}{S}u^{2/\alpha^*}2\Psi\LT(\frac{2u^
{\frac{1}{2}}-\mathbb{Q}_3}{\sqrt{4-\delta(\theta)}}\RT)
=o\LT(u^{2/\alpha^*}\pk{Y(0)>u}\RT).
\end{align*}
Finally by \nelem{in1} for $u$ large enough and $\theta$ small enough
\begin{align*}
\mathcal{A}_4(u)
&\leq \sum_{k=0}^{N(\theta)-1}\LT(\sum_{l=0}^{2M(u)}
\sum_{i=2}^{2M(u)}\qu{J^{k}_{l}(u),J^{k}_{l+i}(u)}
\RT)\\
&\leq\sum_{k=0}^{N(\theta)-1} \sum_{l=0}^{2M(u)}\pk{Y(0)>u}
\LT(\sum_{i=1}^{\IF}\mathbb{Q}_4\exp\LT(-{\mathbb{Q}_5}\abs{iS}^\alpha\RT)\RT)\\
&\leq\mathbb{Q}_6\frac{T}{S}u^{-2/\alpha^*}\pk{Y(0)>u}\LT(\sum_{i=1}
^{\IF}\exp\LT(-{\mathbb{Q}_5}\abs{iS}^\alpha\RT)\RT)\\
&=o\LT(u^{2/\alpha^*}\pk{Y(0)>u}\RT), \ u\rw\IF, S\rw\IF, \theta\rw 0.
\end{align*}
Thus $$\sum_{i=1}^{4}\mathcal{A}_i(u)=o\LT(u^{2/\alpha^*}\pk{Y(0)>u}\RT), \ u\rw\IF,\ S\rw\IF, \theta\rw 0,$$
which derives the result.

\underline{(2)}
In the proof of \netheo{Thm3}, if we take $\beta^*=\frac{2\gamma p}{2-p}$ and  $f(t)=\frac{c t^\gamma}{p}$, then all argumentations in the proof still holds and the results follow.
\COM{
\underline{(3)}\emph{Upper bound}. Set $T_0=0, T_j=\frac{t_j+t_{j+1}}{2},j=1,2,\ldots,n-1,$ and $T_n=T$, then $t_j\in(T_{j-1},T_{j}),j=1,2,\ldots,n$.
By \underline{(1)}, we know
\BQNY
&&\pk{\sup_{t\in[T_{j-1},T_j]}(g(\vk{X}(t))+h(t))>u}\\
&&\quad\quad\sim\pk{g(\vk{X}(0))>u-h_m}
\LT\{
\begin{array}{ll}
2 u^{\frac{2}{\alpha^*}-\frac{2}{\beta^*}}
a_j^{\frac{1}{\alpha}}\mathcal{H}_{\alpha}g^{-\frac{2}{\alpha^*}}\int_{0}^\IF \exp\LT(-\frac{c x^\gamma}{p \gg^{2/p}}\RT)dx,& \text{if}\ \beta^*>\alpha^*,\\
\mathcal{P}_{\alpha,\frac{a_j}{\gg^{2/p}}}^{f}(-\IF,\IF), & \text{if}\ \beta^*=\alpha^*,\\
1,&\text{if}\ \beta^*<\alpha^*.
\end{array}
\RT.
\EQNY
By the Bonferroni inequality we have
\BQNY
\pk{\sup_{t\in[0,T]}(g(\vk{X}(t))+h(t))>u}\leq\sum_{j=1}^n
\pk{\sup_{t\in[T_{j-1},T_j]}(g(\vk{X}(t))+h(t))>u}.
\EQNY
Thus the upper bound follows.\\
\emph{Lower bound}. Since $r(s,t)<1$ for all $s,t\in[0,T],t\neq s$, then
$$\theta:=1-\max_{1\le j, k\le n, j\neq k}r(t_j,t_k)>0.$$
Further, since $r(s,t)$ is continuous function for $(s,t)\in[0,T]\times[0,T]$, we can find small enough $\vn_j, j=1,2,\ldots,n$ such that
$$r(s,t)<1-\frac{\theta}{2},\ \forall\ (s,t)\in[t_j-\vn_j,t_j+\vn_j]\times[t_k-\vn_k,t_k+\vn_k],$$
and $$[t_j-\vn_j,t_j+\vn_j]\cap[t_k-\vn_k,t_k+\vn_k]=\emptyset,$$
for $ 1\leq j,k\leq n, j\neq k.$ Then
\BQNY
&&\E{X_i(s)X_i(t)}=r(s,t)\leq 1-\frac{\theta}{2}, 1\leq i\leq d,
\EQNY
for $(s,t)\in[t_j-\vn_j,t_j+\vn_j]\times[t_k-\vn_k,t_k+\vn_k],1\leq j,k\leq n, j\neq k$.
Set $h_m=\sup_{t\in[0,T]}h_{t}$ and by \nelem{tt},
\BQNY
&&\sum_{1\leq j<k\leq n}\pk{\sup_{t\in[t_j-\vn_j,t_j+\vn_j]}(g(\vk{X}(t))+h(t))>u,
\sup_{t\in[t_k-\vn_k,t_k+\vn_k]}(g(\vk{X}(t))+h(t))>u}\\
&\leq&\sum_{1\leq j<k\leq n}\pk{\sup_{t\in[t_j-\vn_j,t_j+\vn_j]}g(\vk{X}(t))>u-h_m,
\sup_{t\in[t_k-\vn_k,t_k+\vn_k]}g(\vk{X}(t))>u-h_m}\\
&=&o\LT(\pk{g(\vk{X}(0))>u-h_m}\RT).
\EQNY
Using again the Bonferroni inequality we have
\BQNY
\pk{\sup_{t\in[0,T]}(g(\vk{X}(t))+h(t))>u}&\geq&\sum_{j=1}^n
\pk{\sup_{t\in[t_j-\vn_j,t_j+\vn_j]}(g(\vk{X}(t))+h(t))>u}\\
&&-\sum_{1\leq j<k\leq n}\pk{\sup_{t\in[t_j-\vn_j,t_j+\vn_j]}(g(\vk{X}(t))+h(t))>u,
\sup_{t\in[t_k-\vn_k,t_k+\vn_k]}(g(\vk{X}(t))+h(t))>u}.
\EQNY
And we know
\BQNY
\pk{\sup_{t\in[t_j-\vn_j,t_j+\vn_j]}(g(\vk{X}(t))+h(t))>u}
\geq\pk{g(\vk{X}(0))>u-h_m}, \EQNY
so
\BQNY
\pk{\sup_{t\in[0,T]}(g(\vk{X}(t))+h(t))>u}&\geq&
\sum_{j=1}^n\pk{\sup_{t\in[t_j-\vn_j,t_j+\vn_j]}
(g(\vk{X}(t))+h(t))>u}(1+o(1)),\ u\rw\IF.
\EQNY
Combining the above inequality with
\BQNY
&&\pk{\sup_{t\in[t_j-\vn_j,t_j+\vn_j]}(g(\vk{X}(t))+h(t))>u}\\
&&\quad\quad\sim\pk{g(\vk{X}(0))>u-h_m}
\LT\{
\begin{array}{ll}
2 u^{\frac{2}{\alpha^*}-\frac{2}{\beta^*}}
a_j^{\frac{1}{\alpha}}\mathcal{H}_{\alpha}
g^{-\frac{2}{\alpha^*}}\int_{0}^\IF \exp\LT(-\frac{c x^\gamma}{p \gg^{2/p}}\RT)dx,& \text{if}\ \beta^*>\alpha^*,\\
\mathcal{P}_{\alpha,\frac{a_j}{\gg^{2/p}}}^{f}(-\IF,\IF), & \text{if}\ \beta^*=\alpha^*,\\
1,&\text{if}\ \beta^*<\alpha^*,
\end{array}
\RT.
\EQNY
the results follow.\\}

\underline{(3)} If $p=2$, for any constant $\theta>0$, we define
$$I_k=[k\theta, (k+1)\theta],\quad a_k=a(k\theta),\quad  k\in\N, \  N(\theta)= \LT\lfloor\frac{T}{\theta}\RT\rfloor,$$
and $$M^1_\theta(k)=\sup_{t\in I_k}h(t),\ \ \  M^2_\theta(k)=\inf_{t\in I_k}h(t). $$
Then we have
\BQNY
\pu{[0,T]}\geq\sum_{k=0}^{N(\theta)-1}\pu{I_k}-\sum_{j=1}^2\Lambda_j,
\EQNY
where
\BQNY
\Lambda_1=\sum_{k=0}^{N(\theta)}\pu{I_k ,I_{k+1}},\quad
\Lambda_2=\underset{j>k+1}{\sum_{k=0}^{N(\theta)}}\pu{I_k,I_{j}},
\EQNY
and by \underline{(1)}
\begin{align*}
\pu{[0,T]}&\leq\sum_{k=0}^{N(\theta)}\pu{I_k}
\leq\sum_{k=0}^{N(\theta)}\pk{\sup_{t\in I_k}Y(t)>u-M^1_\theta(k)}\\
&\sim \sum_{k=0}^{N(\theta)}a_k^{\frac{1}{\alpha}}\LT(u-M^1_\theta(k)\RT)
^{\frac{1}{\alpha}}\mathcal{H}_\alpha \theta\pk{Y(0)>u-M^1_\theta(k)}\\
&\sim u^{\frac{1}{\alpha}}\mathcal{H}_\alpha \pk{Y(0)>u}\theta \sum_{k=0}^{N(\theta)}a_k^{\frac{1}{\alpha}}e^{\frac{M^1_\theta(k)}{2}}\\
&\sim u^{\frac{1}{\alpha}}\mathcal{H}_\alpha \pk{Y(0)>u}\int_{0}^T (a(t))^{\frac{1}{\alpha}}e^{\frac{h(t)}{2}}d t, \ u\rw\IF,\ \theta\rw 0.
\end{align*}
Similarly,
\begin{align*}
\sum_{k=0}^{N(\theta)-1}\pu{I_k}
&\geq\sum_{k=0}^{N(\theta)-1}\pk{\sup_{t\in I_k}Y(t)>u-M^2_\theta(k)}\\
&\sim u^{\frac{1}{\alpha}}\mathcal{H}_\alpha \pk{Y(0)>u}\int_{0}^T (a(t))^{1/\alpha}e^{\frac{h(t)}{2}}d t, \ u\rw\IF,\ \theta\rw 0.
\end{align*}
Further, we have
\begin{align*}
\Lambda_1&\leq\sum_{k=0}^{N(\theta)}\left(\pu{I_{k}}+\pu{I_{k+1}}-\pu{I_{k}\cup I_{k+1}}\right)\nonumber\\
&\leq\sum_{k=0}^{N(\theta)}\left(\pk{ \sup_{t\in I_{k}}  Y(t)>u-\widetilde{M}^1_\theta(k)}+\pk{\sup_{t\in I_{k+1}} Y(t)>u-\widetilde{M}^1_\theta(k)}-\pk{\sup_{t\in I_{k}\cup I_{k+1}}  Y(t)>u-\widetilde{M}^1_\theta(k)}\right)\nonumber\\
&\sim   \sum_{k=0}^{N(\theta)}\LT(a_k^{1/\alpha}+a_{k+1}^{1/\alpha}
-2a_k^{1/\alpha}\RT)\theta
u^{\frac{1}{\alpha }}e^{\frac{\widetilde{M}^1_\theta(k)}{2}}\pk{Y(0)>u}
=o\left(u^{1/\alpha}\pk{Y(0)>u}\right), \ u\rw\IF,\ \theta\rw 0,
\end{align*}
where $\widetilde{M}^1_\theta(k)=\max(M^1_\theta(k),M^1_\theta(k+1))$.\\
Set $h_m=\sup_{t\in[0,T]}h(t)$ and for any $\theta>0$
\BQNY
\E{X_i(t)X_i(s)}=r(|t-s|)\leq 1-\delta(\theta)
\EQNY
for $(s,t)\in I_k\times I_j,j>k+1$ where $\delta(\theta)>0$ is a constant related to $\theta$. Then by \nelem{tt}
\begin{align*}
\Lambda_2&=\underset{j>k+1}{\sum_{k=0}^{N(\theta)}}\pu{I_k,I_{j}}\leq \underset{j>k+1}{\sum_{k=0}^{N(\theta)}}\pk{\sup_{t\in I_k}Y(t)>u-h_m,\sup_{t\in I_{j}}Y(t)>u-h_m}\\
&\leq\underset{j>k+1}{\sum_{k=0}^{N(\theta)}}2\Psi\LT(\frac{2\LT(u-h_m\RT)
^{\frac{1}{2}}-\mathbb{Q}_1}{\sqrt{4-\delta(\theta)}}\RT)=o\LT(\pk{Y(0)>u}\RT), \ u\rw\IF, \theta\rw0.
\end{align*}
Thus, we have
\BQNY
\pu{[0,T]}\sim \int_{0}^T (a(t))^{1/\alpha}e^{\frac{h(t)}{2}}d t\mathcal{H}_\alpha u^{\frac{1}{\alpha}} \pk{Y(0)>u}, \ u\rw\IF.
\EQNY
If $p\in (2,\IF)$, set $M_1=\inf_{t\in[0,T]}h(t)$ and $M_2=\sup_{t\in[0,T]}h(t)$. Since $h(t)$ is a continuous function, we have
$-\IF<M_1\leq M_2<\IF$.
Further, since when $p\in(2,\IF)$, $$\pk{Y(0)>u+\mathbb{Q}_2}\sim\pk{Y(0)>u}$$
hold for any $\mathbb{Q}_2\in\R$.
Hence, by \underline{(1)}
\BQNY
\pu{[0,T]}\geq \pk{\sup_{t\in[0,T]}Y(t)>u-M_1}
\sim \int_{0}^T (a(t))^{1/\alpha}dtu^{\frac{2}{\alpha p}}\mathcal{H}_\alpha \pk{Y(0)>u}, \ u\rw\IF,
\EQNY
and
\BQNY
\pu{[0,T]}\leq \pk{\sup_{t\in[0,T]}Y(t)>u-M_2}
\sim\int_{0}^T (a(t))^{1/\alpha}dtu^{\frac{2}{\alpha p}}\mathcal{H}_\alpha \pk{Y(0)>u}, \ u\rw\IF.
\EQNY
Thus we complete the proof.

\QED

\COM{
\section{Appendix: Uniform Pickands-Piterbarg theorem on short intervals}
In this section we are interested on the Pickands-Piterbarg Theorem which establishes the asymptotics of supremum of non-Gaussian random process in a short interval under regularity conditions on the correlation function.
For a compact set $E\in \R$, let $\{g_u(t),t\in E\}$ be a family of general random processes with continuous trajectories.\\
Next, let $\Gamma: C(E)\to \R$   be a real-valued functional such that for all $u$ large:
As in \cite{dekos14} we shall assume: \\
{\bf F1:} $\Gamma(af+b)=a\Gamma(f)+b$ for any $f\in C_0(E)$ and $a>0,b\inr$ and further $\Gamma$ is continuous.\\
{\bf F2:} There exits $c>0$ such that  $\abs{\Gamma(f)}\le c \sup_{{t}\in E} {\abs{f({t})}}$ for any $f\in C_0(E)$. \\
Further, we assume that $K_u$ a family of index sets and $u_k$ satisfying that
\BQN\label{uk1}
\lim_{u\rw\IF}\sup_{k\in K_u}\LT|\frac{u_k}{u}-1\RT|=0.
\EQN

For all $u$ large enough
\BQN\label{xhoP}
\frac{ \pk{\Gamma(g_u(t))> u_k } } {  \pk{ V> u_k}}
&=&\int_{\R} \pk{\Gamma(g_u(t))> u_k \Bigl \lvert
g_u(0) = u_k-d(u_k)w}d H_{V,u_k}(u_k-d(u_k)w)\nonumber\\
&=&\int_{\R} \pk{Z_{u_k}(w)>0}d H_{V,u_k}(u_k-d(u_k)w)
\EQN
where  $H_{V,u_k}(u_k-d(u_k)w)= \frac{1-F_V(u_k-d(u_k)w)}{1-F_V(u_k)}$ and $Z_{u_k}(w):=\Gamma(g_u(t))-u_k \Bigl \lvert
g_u(0) = u_k-d(u_k)w$.

We shall impose the following assumptions:\\
 {\bf A1:} $V:= g_u(0)$ does not depend on $u$ with distribution function $F_V$ and density function $f_V$ which
satisfies
\BQNY
\lim_{u\rw\IF}\sup_{k\in K_u, w<W}\LT|\frac{d(u_k)f_V(u_k-d(u_k)w)}{(1-F_V(u_k))e^{w}}-1\RT|=0,
\EQNY
where  $W>0$ is a large constant and $d(x)$ is a function with $\lim_{x\rw\IF}d(x)/x=0$.\\
{\bf A2:}
There exists a random process $\{\Theta(t), t\in E\}$ with
\BQN\label{gaupper}
\lim_{W\rw \IF}\int_{W}^{\IF}e^{w}\pk{\Gamma(\Theta(t))>w}d w=0
\EQN
such that for any large $W>0$
\BQN \label{convergence}
\lim_{u\rw\IF}\sup_{k\in K_u, w\in[-W,W]}\LT| \pk{Z_{u_k}(w)>0}-\pk{\Gamma(\Theta(t))>w}\RT|=0.
\EQN
{\bf A3:} For some $u_0>0$
\BQN
\limit{W} \sup_{u>u_0}\sup_{k\in K_u}\int_{W}^\IF\pk{Z_{u_k}(w)>0} d H_{V,u_k}(u_k-d(u_k)w)=0.
\EQN

\BT \label{th2}
For $g_u(t), \Gamma, u_k$ above satisfy assumptions {\bf A1}--{\bf A3}, we have
\BQN
\limit{u} \sup_{k\in K_u}\LT|\frac{ \pk{\Gamma(g_u(t))> u_k } } {  \pk{ V> u_k}} -\int_\R e^w \pk{\Gamma(\Theta(t))> w} d w\RT|=0.
\EQN
\ET

\BRM
i) If $g_u(0) \le \Gamma(g_u(t))$, which is in particular the case if
$\Gamma(f)=\sup_{t\in E} f(t)$, then
we have that
$$ \pk{\Gamma(g_u(t))> u}= \pk{g_u(0)> u}+ \pk{\Gamma(g_u(t))> u, g_u(0) \le u}
$$
implying that
$$\int_\R e^w \pk{\Gamma(\Theta(t))> w} \, dw = 1+ \int_0^\IF  e^w \pk{ \Gamma(\Theta(t))> w} d w.$$

\ERM
\prooftheo{th2}
 $$\frac{d H_{V,u_k}(u_k-d(u_k)w)}{d w}= \frac{d(u_k)f_V(u_k-d(u_k)w)}{1-F_V(u_k)}.$$
For all $u$ large enough and any $W>0$
\BQN\label{xhoP}
\frac{ \pk{\Gamma(g_u(t))> u_k } } {  \pk{ V> u_k}}
&=&\int_{\R} \pk{\Gamma(g_u(t))> u_k \Bigl \lvert
g_u(0) = u_k-d(u_k)w}d H_{V,u_k}(u_k-d(u_k)w)\nonumber\\
&=&\int_{-\IF}^{-W} \pk{Z_{u_k}(w)>0}d H_{V,u_k}(u_k-d(u_k)w)+\int_{-W}^{W} \pk{Z_{u_k}(w)>0}d H_{V,u_k}(u_k-d(u_k)w)\nonumber\\
&& +\int_{W}^{\IF} \pk{Z_{u_k}(w)>0}d H_{V,u_k}(u_k-d(u_k)w)\nonumber\\
&=:&I_1(u)+I_2(u)+I_3(u),
\EQN
and
\BQNY
&&\sup_{k\in K_u}\LT|\frac{ \pk{\Gamma(g_u(t))> u_k } } {  \pk{ V> u_k}} -\int_\R e^w
\pk{ \Gamma(\Theta(t))> w} d w\RT|\\
&&\quad\quad\leq \sup_{k\in K_u}\LT|I_1(u)\RT|+\sup_{k\in K_u}\LT|I_2(u) -\int_{-W}^{W} e^w \pk{\Gamma(\Theta(t))> w} d w\RT|+\sup_{k\in K_u}\LT|I_3(u)\RT|\\
&&\quad\quad\quad+\int_{-\IF}^{-W} e^w \pk{\Gamma(\Theta(t))> w} d w+\int_{W}^{\IF} e^w \pk{\Gamma(\Theta(t))> w} d w
\EQNY
By assumption {\bf A1}
\BQNY
\sup_{k\in K_u}\LT|I_1(u)\RT|&=&\sup_{k\in K_u}\int_{-\IF}^{-W} \pk{Z_{u_k}(w)>0}d H_{V,u_k}(u_k-d(u_k)w) \\
&=&\sup_{k\in K_u}\int_{-\IF}^{-W} \pk{Z_{u_k}(w)>0}\frac{d(u_k)f_V(u_k-d(u_k)w)}{(1-F_V(u_k))} d w \\
&\leq& \mathbb{Q}\int_{-\IF}^{-W} e^{w}d w\\
& \le& \mathbb{Q}e^{-W}\rw 0, \quad u \to \IF , W\rw\IF,
\EQNY
and by assumption {\bf A3}
\BQNY
\sup_{k\in K_u}\LT|I_3(u)\RT|=\sup_{k\in K_u}\int_{W}^{\IF} \pk{Z_{u_k}(w)>0}d H_{V,u_k}(u_k-d(u_k)w) \rw 0, \quad u \to \IF , W\rw\IF.
\EQNY
Furthermore, by assumption {\bf A1}
\BQNY
\sup_{k\in K_u, w\in[-W,W]}\LT|\frac{d(u_k)f_V(u_k-d(u_k)w)}{(1-F_V(u_k))}-e^{w}\RT|&\leq& e^W\sup_{k\in K_u, w\in[-W,W]}\LT|\frac{d(u_k)f_V(u_k-d(u_k)w)}{(1-F_V(u_k))e^{w}}-1\RT|\\
&\rw& 0, \ u\rw\IF,
\EQNY
which combined with  assumption {\bf A2} follows by
\BQNY
\lim_{u\rw 0}\sup_{k\in K_u}\LT|I_2(u) -\int_{-W}^{W} e^w \pk{\Gamma( \Theta(t))> w} d w\RT|=0.
\EQNY
By \eqref{gaupper}
\BQNY
\int_{-\IF}^{-W} e^w \pk{\Gamma(\Theta(t))> w} d w+\int_{W}^{\IF} e^w
\pk{\Gamma(\Theta(t))> w} d w\leq e^{-W} +\int_{W}^{\IF} e^w
\pk{\Gamma(\Theta(t))> w} d w\rw 0,\ W\rw\IF,
\EQNY
establishing the proof.
\QED}
\section{Appendix}
\subsection{Appendix A}
First we give the proof of \eqref{intervalmax}.

{\bf Proof of \eqref{intervalmax}:}
\cLb{Here we use $\pu{\cdot}$ the same as in \eqref{ppu}.}  We consider the case $0<A<B<T$. First by \underline{(1)} of \netheo{Thm1}, we have
 as $u\rw\IF$
\BQNY
\pu{[A,B]}=\pk{\sup_{t\in[A,B]} Y(t)>u-h_m}
\sim\int_{A}^B(a(t))^{\frac{1}{\alpha}}dt\mathcal{H}_\alpha u^{\frac{2}{\alpha^*}}\pk{Y(0)>u-h_m}.
\EQNY
Set $\Delta_\vn=[A-\vn, B+\vn]\cap [0,T]$ for some $\vn>0$, then we have
\BQNY
\pu{[A,B]}\leq\pu{[0,T]}\leq \pu{[\Delta_\vn]}+\pu{[0,T]\setminus\Delta_\vn}.
\EQNY
 Since $h(t)$ is a continuous function and we have $h_\vn:=\sup_{t\in [0,T]\setminus\Delta_\vn}h(t)<h_m$, then by \underline{(1)} of \netheo{Thm1}
\begin{align*}
\pu{[0,T]\setminus\Delta_\vn}&\leq \pk{\sup_{t\in[0,T]\setminus\Delta_\vn} Y(t) >u-h_\vn}
\sim \mathbb{Q}_1u^{2/\alpha^*} \pk{Y(0)>u-h_\vn}\\
&=o\LT(u^{2/\alpha^*}\pk{Y(0)>u-h_m}\RT),\ u\rw\IF,\ \vn\rw 0.
\end{align*}
Further, we have by \underline{(1)} of \netheo{Thm1}
\begin{align*}
\pu{[\Delta_\vn]}&\leq\pk{\sup_{t\in\Delta_\vn} Y(t)>u-h_m}
\leq\int_{A-\vn}^{B+\vn}(a(t))^{\frac{1}{\alpha}}dt\mathcal{H}_\alpha u^{\frac{2}{\alpha^*}}\pk{Y(0)>u-h_m}\\
&\sim\int_{A}^B(a(t))^{\frac{1}{\alpha}}dt\mathcal{H}_\alpha u^{\frac{2}{\alpha^*}}\pk{Y(0)>u-h_m},\  u\rw\IF,\ \vn\rw 0.
\end{align*}
Hence the claims follow.

\QED

{\bf proof of Example \ref{ex1}:} Notice that we just need to prove that $g(\vk{x})$ satisfies \necon{gcon}.\\
\COM{If $\rho=2$, then we have $\gg=1$ and $\mathcal{M}$ is the whole unit sphere $\mathbb{S}_{d-1}$ (which is a manifold of dimension $d-1$). According to \cite{ELPNT2017}, we know that \netheo{Thm3} and \netheo{Thm1} hold.\\
If $\rho<2$, then $$}
i) By the twice continuously differentiable of $g$, it easily to prove that \necon{gcon} is satisfied, except case $\rho=2$. For the case $\rho=2$, we have $\gg=1$ and $\mathcal{M}$ is the whole unit sphere $\mathbb{S}_{d-1}$ (which is a manifold of dimension $d-1$), and according to \cite{ELPNT2017}, we know that \netheo{Thm3} and \netheo{Thm1} hold. A similar detail analysis can be found in \cite{chaos15}[Example 1].\\
\COM{For $\rho=2$, we have $\gg=1$ and $\mathcal{M}=\{ \vk{v}: (v_1\ldot v_m)\in \mathcal{S}_{m-1}, v_i=0, m+1\leq i\leq d\}$ (which is a manifold of dimension $m-1$). We introduce the spherical coordinates
\BQNY
&& v_1= \cos \varphi_1,\\
&& v_2= \sin \varphi_1\cos \varphi_2,\\
&&\ldots\\
&& v_{d-1}= \sin \varphi_1\sin \varphi_2\ldots \sin \varphi_{d-2}\cos \varphi_{d-1},\\
&& v_{d}= \sin \varphi_1\sin \varphi_2\ldots \sin \varphi_{d-2}\sin \varphi_{d-1},
\EQNY
with $\vk{\varphi}\in\Pi_{d-1}=[0,\pi)^{d-2}\times[0,2\pi)$. Then for any $\vk{\varphi}\in\mathcal{M}$, we have}
ii) For $g(\vk{x})=\Pi_{i=1}^dx_i$, we have $p=d$ and further $\gg=\frac{1}{d^{d/2}}$ since
\BQNY
\mathcal{M}=\LT\{\LT(\pm1/\sqrt{d}\ldot \pm1/\sqrt{d} \RT) \text{with even number of negative coordinates} \RT\},
\EQNY
which consists of $2^{d-1}$ points (the product $x_1\ldot x_d$ should be  positive). Further, we introduce the spherical coordinates
\BQNY
&& v_1= \cos \varphi_1,\\
&& v_2= \sin \varphi_1\cos \varphi_2,\\
&&\ldots\\
&& v_{d-1}= \sin \varphi_1\sin \varphi_2\ldots \sin \varphi_{d-2}\cos \varphi_{d-1},\\
&& v_{d}= \sin \varphi_1\sin \varphi_2\ldots \sin \varphi_{d-2}\sin \varphi_{d-1},
\EQNY
with $\vk{\varphi}\in\Pi_{d-1}=[0,\pi)^{d-2}\times[0,2\pi)$ and
\BQNY
g(\vk{\varphi})=\sin^{d-1}\varphi_1\ldot \sin\varphi_{d-1} \cos \varphi_1\ldot  \cos \varphi_{d-1}.
\EQNY
For instance, at the point $\LT(1/\sqrt{d}\ldot 1/\sqrt{d} \RT)$ we have $\cos \varphi_i=\sqrt{\frac{1}{d-i+1}}$ and $\sin \varphi_i=\sqrt{\frac{d-i}{d-i+1}}$. Further calculation shows that for any $\vk{\varphi}\in \mathcal{M}_{\varphi}$
\BQNY
g''_{\varphi_i\varphi_i}(\vk{\varphi})=-2g(\vk{\varphi})(d-i+1)
=-\frac{2(d-i+1)}{d^{d/2}},\quad g''_{\varphi_i\varphi_j}(\vk{\varphi})=0,\
\text{for}\  i\neq j,
\EQNY
which leads $\abs{det g''(\vk{\varphi})}=2^{d-1}d!/d^{d(d-1)/2}>0$.
By the continuity of $\abs{det g''(\vk{\varphi})}$ over some $\mathcal{M}_{\varphi}(\vn)$, we know that \necon{gcon} is satisfied.\\
iii) For $g(\vk{x})=\max_{1\leq i\leq d}x_i$, we have $p=1$ and further $\gg=1$ since
$\mathcal{M} $ consists of $d$ points $\LT(0\ldot 0, 1, 0\ldot 0\RT)$.
For instance, around the point $\vk{v}=\LT(0\ldot 0, 1\RT)$ over unit sphere $\mathbb{S}_{d-1}$ we have that $$g(\vk{\varphi})=\sin \varphi_1\sin \varphi_2\ldots \sin \varphi_{d-2}\sin \varphi_{d-1}$$ and $ \varphi_i=\frac{\pi}{2}, i=1\ldot d$. It is clear that
\BQNY
g''_{\varphi_i\varphi_i}(\vk{\varphi})=-g(\vk{\varphi})
=-1,\quad g''_{\varphi_i\varphi_j}(\vk{\varphi})=0,\
\text{for}\  i\neq j,
\EQNY
which leads $\abs{det g''(\vk{\varphi})}=1>0$.
By the continuity of $\abs{det g''(\vk{\varphi})}$ over some $\mathcal{M}_{\varphi}(\vn)$, we know that \necon{gcon} is satisfied.
\QED

Next we give the proofs of \nelem{tt}, \nelem{im1} and \nelem{in1}.  In the following proofs, $\mathbb{Q}_i,i\in\N$ denote some positive constants which can be different from line by line.

\prooflem{tt}
By \eqref{lrr21} and the continuity of $r(s,t)$, for some $\delta>0$ we have
\BQNY
\E{X_i(t)X_i(s)}=r(s,t)\leq1-\frac{\delta}{2}, i=1,2,\ldots,d,
\EQNY
holds for any $(s,t)\in[0,t_1]\times[t_2,t_3]$.
Set $Z(t,\vk{v},s, \vk{w})=\langle\vk{X}(t),
\vk{v}\rangle+\langle\vk{X}(s),\vk{w}\rangle$ where $\vk{v},\vk{w} \in\mathbb{S}_{d-1}$  with $\mathbb{S}_{d-1}$ the unit sphere in $\R^d$ with respect to Euclidean norm $\norm{\cdot}$. 
Since  $Z(t,\vk{v},s, \vk{w})$ is a center Gaussian fields, we have further
\begin{align*}
\Var\LT(Z(t,\vk{v},s, \vk{w})\RT)&=2+2r(s,t)\LT(\sum_{i=1}^d v_i w_i\RT)
\leq2+r(s,t)\LT(\sum_{i=1}^d v_i^2+ \sum_{i=1}^d w^2_i\RT)\\
&=2+2r(s,t)\leq4-\delta
\end{align*}
for any $(t,\vk{v},s,\vk{w})\in[0,t_1]\times\mathbb{S}_{d-1}\times
[t_2,t_3]\times\mathbb{S}_{d-1}$.
By Borell inequality (see e.g., \cite{AdlerTaylor, GennaBorell})
\begin{align*}
\pk{\sup_{t\in[0,t_1]}Y(t)>u,\sup_{t\in[t_2,t_3]}Y(t)>u}
&\leq\pk{\sup_{t\in[0,t_1]}|Y(t)|>u,
\sup_{t\in[t_2,t_3]}|Y(t)|>u}\\
&\leq\pk{\sup_{t\in[0,t_1]}|\vk{X}(t)|>u^{},
\sup_{t\in[t_2,t_3]}|\vk{X}(t)|>u^{}}\\
&\leq\pk{\sup_{(t,\vk{v})\in[0,t_1]
\times\mathbb{S}^1_{d-1}}\langle\vk{X}(t),
\vk{v}\rangle>u^{},\sup_{(s,\vk{w})\in[t_2,t_3]
\times\mathbb{S}^2_{d-1}}\langle\vk{X}(s),\vk{w}\rangle>
u^{}}\\
&\leq\pk{\sup_{(t,\vk{v},s,\vk{w})\in[0,t_1]
\times\mathbb{S}^1_{d-1}
\times[t_2,t_3]\times\mathbb{S}^2_{d-1}}Z(t,\vk{v},s, \vk{w})>2u^{}}\\
&\leq 2\Psi\LT(\frac{2u^{}-D}{\sqrt{4-\delta}}\RT),
\end{align*}
where $D$ is some constant such that
$$\pk{\sup_{(t,\vk{v},s,\vk{w})\in[0,t_1]\times\mathbb{S}^1_{d-1}\times[t_2,t_3]
\times\mathbb{S}^2_{d-1}}Z(t,\vk{v},s, \vk{w})>D}\leq \frac{1}{2},$$
 hence the first claim follows. The second claim follows for \eqref{tdf}.
\QED

Before giving the proof of \nelem{im1}, we remark that from \cite{PiterChaos2015}[Corollary 4] it follows that only arbitrary small vicinity in the sphere of the maximum point set gives contribution to the asymptotic behavior of probabilities of \eqref{e1}. Therefore we can change $g(\vk{\varphi})$ outside of $\mathcal{M}_{\vk{\varphi}}(\vn)$ is such a way that first,
$$g(\vk{\varphi})\geq 1-\vn, \vk{\varphi}\in\Pi_{d-1},$$
keeping the same asymptotic behavior of the probability in question; and second, having $g(\vk{\varphi})$ twice continuously differentiable for all  $\vk{\varphi}\in\Pi_{d-1}$.\\

\prooflem{im1}
i) First we prove \eqref{pp1} and assume for simplicity that $f_i=f, i\le d$.  For any $W>0$ and all $u$ large (set  $Z_u(t) =\frac{g(\vk{\xi}
	(u^{-2/\alpha}t+t_0))}{1+u^{-2}f(t)}$  and write simple $p(x)$ instead of $\ehe{p}(x)$)
\BQNY
&&\pk{\sup_{t\in[-S_1,S_2]} \ehe{Z_u(t)} >u_k}\\
&&\quad\quad=\int_{-\IF}^{\IF}\pk{\sup_{t\in[-S_1,S_2]}
\ehe{Z_u(t)} >u_k\Big|g(\vk{\xi}(t_0))=x}
\ehe{p}(x)dx\\
&&\quad\quad=u_k^{-1}\int_{-\IF}^{\IF}\pk{\sup_{t\in[-S_1,S_2]}
\ehe{Z_u(t)}>u_k
\Big|g(\vk{\xi}(t_0))=\uy}
\ehe{p}(\uy)dy\\
&&\quad\quad=u_k^{-1}\int_{-W}^{W}\pk{\sup_{t\in[-S_1,S_2]}
\ehe{Z_u(t)} >u_k\Big|g(\vk{\xi}(t_0))=\uy}
\ehe{p}(\uy)dy\\
&&\quad\quad\quad+u_k^{-1}\int_{W}^{\IF}\pk{\sup_{t\in[-S_1,S_2]}
\ehe{Z_u(t)}>u_k\Big|g(\vk{\xi}(t_0))=\uy}
\ehe{p}(\uy)dy\\
&&\quad\quad\quad+\pk{g(\vk{\xi}(t_0))> u_k+u_k^{-1}W}\\
&&\quad\quad=: I_1(u)+I_2(u)+I_3(u),
\EQNY
where \cLa{$\uy=u_k-u_k^{-1}y$} and $\ehe{p}(x)$ is the density function of $g(\vk{\xi}(t_0))$ which is showed in \netheo{Them00}.\\
Since
\BQNY
\mathcal{P}_{\alpha,a}^{f(t)}
[-S_1,S_2]=\E{\exp\LT(\sup_{t\in[-S_1,S_2]}\eta(t)\RT)}
=\int_{-\IF}^{\IF}e^{y}\pk{\sup_{t\in[-S_1,S_2]}\eta(t)>y} d y,
\EQNY
where $\eta(t):=\sqrt{2a}B_{\alpha}(t)-a\abs{t}^\alpha-f(t)$, we have
\BQNY
&&\sup_{k\in K_u}\LT|\frac{\pk{\sup_{t\in [-S_1,S_2]}\ehe{Z_u(t)}>u_k}}
{\pk{g(\vk{\xi}(t_0))>u_k}}-\mathcal{P}_{\alpha,a}
^{f(t)}
[-S_1,S_2]\RT|\\
&&\leq\sup_{k\in K_u}\LT|\frac{I_1(u)}
{\pk{g(\vk{\xi}(t_0))>u_k}}-\mathcal{P}_{\alpha,a}
^{f(t)}
[-S_1,S_2]\RT|+\sup_{k\in K_u}\frac{I_2(u)}
{\pk{g(\vk{\xi}(t_0))>u_k}}+\sup_{k\in K_u}\frac{I_3(u)}
{\pk{g(\vk{\xi}(t_0))>u_k}}\\
&&\leq\sup_{k\in K_u}\LT|\frac{I_1(u)}
{\pk{g(\vk{\xi}(t_0))>u_k}}-\int_{-W}^{W}e^{y}\pk{\sup_{t\in[-S_1,S_2]}\eta(t)>y} d y\RT|\\
&&\quad+\int_{\abs{y}\geq W}e^{y}
\pk{\sup_{t\in[-S_1,S_2]}\eta(t)>y} d y
+\sup_{k\in K_u}\frac{I_2(u)}
{\pk{g(\vk{\xi}(t_0))>u_k}}+\sup_{k\in K_u}\frac{I_3(u)}
{\pk{g(\vk{\xi}(t_0))>u_k}}.
\EQNY
By \cite{Pit96} [Theorem 8.1] and the fact that for $s,t\geq 0$
$$
\E{\LT(B_\alpha(t)-B_\alpha(s)\RT)^2}=\abs{t-s}^\alpha,$$
we have for any $y>0$
\begin{align}\label{boundeta}
\pk{\sup_{t\in[-S_1,S_2]}\eta(t)>y}&\leq
\pk{\sup_{t\in[-S_1,S_2]} \sqrt{2a}B_{\alpha}(t)>y}
\leq 2\pk{\sup_{t\in[0,\max(S_1,S_2)]} B_{\alpha}(t)>\frac{y}{\sqrt{2a}}}\nonumber\\
&\leq \mathbb{Q}_0y^{2/\alpha}\Psi\LT(\frac{y}{\sqrt{2a\max(S_1,S_2)}}\RT),
\end{align}
which implies that
\begin{align*}
\int_{\abs{y}>W}e^{y}\pk{\sup_{t\in[-S_1,S_2]}\eta(t)>y} d y
&\leq \int_{-\IF}^{-W}e^{y}\pk{\sup_{t\in[-S_1,S_2]}\eta(t)>y} d y+
\int_{W}^{\IF}e^{y}\pk{\sup_{t\in[-S_1,S_2]}\eta(t)>y}d y\\
&\leq e^{-W}+\int_{W}^{\IF}\mathbb{Q}_0e^{y} y^{2/\alpha}\Psi\LT(\frac{y}{\sqrt{2a\max(S_1,S_2)}}\RT) d y
\rw0, \ W\rw\IF.
\end{align*}
Next, using \eqref{tdf} and \eqref{uk}, we have
\BQN\label{II3}
\sup_{k\in K_u}\frac{I_3(u)}{\pk{g(\vk{\xi}(t_0))>u_k}}\leq \sup_{k\in K_u}\frac{\pk{ g(\vk{\xi}(t_0))> u_k+u_k^{-1}W}}{\pk{g(\vk{\xi}(t_0))>u_k}}\rw 0,\quad u\rw\IF,\quad W\rw \IF,
\EQN
where we use the fact that
$$\LT(u_k+u_k^{-1}W\RT)^{2}=u_k^{2}+2W+u_k^{-2}W^2.$$
By \eqref{uk}, \eqref{pdf} and \eqref{tdf} for $\vn_1<\frac{1}{2}$
\BQN\label{bb}
&&\sup_{k\in K_u, -W<y<u_k^{2-\vn_1}}\LT|\frac{u_k^{-1}
\ehe{p}(\uy)}
{\pk{g(\vk{\xi}(t_0))>u_k}
e^{y}}-1\RT|\nonumber\\
&&=\sup_{k\in K_u, -W<y<u_k^{2-\vn_1}}\LT|\frac{u_k^{-1}
h_0(1+R(\uy))
\LT(\uy\RT)^{m}
e^{-\frac{1}{2}\uy^{2}}}
{\pk{g(\vk{\xi}(t_0))>u_k}
e^{y}}-1\RT|\nonumber\\
&&=\sup_{k\in K_u, -W<y<u_k^{2-\vn_1}}\LT|\frac{(1+R(\uy))\LT(1-u_k^{-2}y\RT)
^{m}e^{-\frac{y^2}{2 u_k^{2}}}}
{1+R_1(u_k)}-1\RT|\rw 0,\  u\rw\IF,
\EQN
where $R(x)\rw 0, R_1(x)\rw 0, x\rw \IF$, $h_0$ is the same as in \netheo{Them00}, and we use the fact that
\BQNY
\LT(\uy\RT)^{2}=u_k^{2}\LT(1-u_k^{-2}y\RT)^{2}=u_k^{2}
-2y+u_k^{-2}y^2,
\EQNY
and
\BQNY
\LT(\uy\RT)^{m}
=u_k^{m}(1-u_k^{-2}y)^{m}.
\EQNY
Setting $$P_{u_k}(y):=\pk{\sup_{t\in[-S_1,S_2]}
\ehe{Z_u(t)}>u_k
\Big|g(\vk{\xi}(t_0))=\uy}\leq 1,$$
we have
\BQN\label{basic}
&&\sup_{k\in K_u}\LT|\frac{I_1(u)}
{\pk{g(\vk{\xi}(t_0))>u_k}}-\int_{-W}^{W}e^{y}
\pk{\sup_{t\in[-S_1,S_2]}\eta(t)>y} d y\RT|\nonumber\\
&&\leq \sup_{k\in K_u}\LT|\int_{-W}^{W}P_{u_k}(y)
\frac{u_k^{-1}\ehe{p}(\uy)}
{\pk{g(\vk{\xi}(t_0))>u_k}}dy-\int_{-W}^{W}e^{y}P_{u_k}(y) d y\RT|\nonumber\\
&&\quad +\sup_{k\in K_u}\int_{-W}^{W}\abs{e^{y}P_{u_k}(y)-e^{y}\pk{\sup_{t\in[-S_1,S_2]}\eta(t)>y}}d y\nonumber\\
&&\leq \sup_{k\in K_u}e^{W}\int_{-W}^{W}\LT|
\frac{u_k^{-1}\ehe{p}(\uy)}
{e^{y}\pk{g(\vk{\xi}(t_0))>u_k}}-1\RT|d y +\sup_{k\in K_u}e^{W}\int_{-W}^{W}\abs{P_{u_k}(y)-\pk{\sup_{t\in[-S_1,S_2]}\eta(t)>y}}d y,
\EQN
where by \eqref{bb}, the first term goes to zero as $u\rw\IF$.\\
Thus we just need to prove
 \BQNY
\lim_{u\rw\IF}\sup_{k\in K_u, y\in[-W,W]}\LT|P_{u_k}(y)-\pk{\sup_{t\in[-S_1,S_2]}\eta(t)>y}\RT|=0
\EQNY
and
\BQNY
\lim_{u\rw\IF}\sup_{k\in K_u}\frac{I_2(u)}
{\pk{g(\vk{\xi}(t_0))>u_k}}=0,
\EQNY
which will be dealt with in the following {\bf Step 1} and {\bf Step 2}, respectively.

{\bf Step 1}:
Using the idea of Piterbarg, we use  $\vk{e}\in \mathbb{S}_{d-1}$ as the spherical coordinates of $\vk{v}$ with
$$\mathcal{V}:=\{\vk{v}:g(\vk{v})=1,\vk{v}\in\R^d\}=\{\vk{v}: \vk{v}=\vk{v}(\vk{e})=|\vk{v}|\vk{e}, \vk{e}\in \mathbb{S}_{d-1}, g(\vk{v})=1\},$$
and write
\BQNY
\{g(\vk{\xi}(t_0))=\uy\}=:\bigcup_{\vk{e}\in \mathbb{S}_{d-1}} U_{\vk{e}},
\EQNY
where $U_{\vk{e}}=\{(\uy)^{-1}\vk{\xi}(t_0)=\vk{v}(\vk{e})\}$.
Then
\BQN\label{pku}
P_{u_k}(y)=\int_{\mathbb{S}_{d-1}}f_k(\vk{e})\pk{\sup_{t\in[-S_1,S_2]}
\ehe{Z_u(t)}>u_k
\Big|U_{\vk{e}}}d\vk{e}
\EQN
where $d\vk{e}$ is the elementary volume in $\mathbb{S}_{d-1}$, $f_k$ is the conditional probability density of the Gaussian vector
$(\uy)^{-1}\vk{\xi}(t_0)$ given $g((\uy)^{-1}\vk{\xi}(t_0))=1$ induced on $\mathbb{S}_{d-1}$.\\
We have
\begin{align*}
\E{\frac{\xi_i(u^{-2/\alpha}t+t_0)}{1+u^{-2}f(t)}\Big|U_{\vk{e}}}
=\E{\frac{\xi_i(u^{-2/\alpha}t+t_0)}{1+u^{-2}f(t)}
\Big|\xi_i(t_0)=\uy v_i}
=\frac{r(u^{-2/\alpha}|t|)\uy v_i}{1+u^{-2}f(t)}
=: m_u^i(t),
\end{align*}
$i=1,\ldots,d$, where $\vk{v}=\vk{v}(\vk{e})=(v_1,\ldots,v_d)$. Further calculating the conditional variance we have
\BQNY
&&\Var\LT(\frac{\xi_i(u^{-2/\alpha}s+t_0)}{1+u^{-2}f(s)}
-\frac{\xi_i(u^{-2/\alpha}t+t_0)}{1+u^{-2}f(t)}\Big|U_{\vk{e}}\RT)\\
&&\quad\quad=\Var\LT(\frac{\xi_i(u^{-2/\alpha}s+t_0)}{1+bu^{-2}f(s)}
-\frac{\xi_i(u^{-2/\alpha}t+t_0)}{1+bu^{-2}f(t)}\RT)-
\LT(\frac{r(u^{-2/\alpha}|s|)}{1+u^{-2}f(s)}
-\frac{r(u^{-2/\alpha}|t|)}{1+u^{-2}f(t)}\RT)^2\\
&&\quad\quad=:V_u(s,t).
\EQNY
Let us recall that
\BQNY
\uy=u_k\LT(1-y u_k^{-2}\RT),
\EQNY
which combined with \eqref{lrr11} leads to
\BQN
&&m_u^i(t)=u_k v_i\LT(1-\LT(f(t)+a|t|^\alpha\RT)u^{-2}
-y u_k^{-2}+u^{-2}\rho_1(u,k,t,y)\RT)\label{mean}\\
&&V_u(s,t)=2au^{-2}(|t-s|^\alpha+\rho_2(u,s,t)),\label{var}
\EQN
where $\rho_1(u,k,t,y)\rw 0$, as $u\rw \IF$ uniformly in $t, y, k$ and $\rho_2(u,s,t)\rw 0$, as $u\rw \IF$ uniformly in $s,t$ .
Let us introduce $d$ independent Gaussian processes $\mathcal{X}^i_u(t), i=1,2,\ldots,d$ with zero means, continuous trajectories with $t\in[-S_1,S_2]$, and covariance functions equal to the respective conditional covariance function of $\frac{\xi_i(u^{-2/\alpha}t+t_0)}{1+u^{-2}f(t)}$. Now write
\BQNY
P_{u_k}(y)=\int_{\mathbb{S}_{d-1}}f_k(\vk{e})\pk{\sup_{t\in[-S_1,S_2]}
g(\vk{\mathcal{X}}_u(t)+\vk{m}_u(t))>u_k}d\vk{e},
\EQNY
with $\vk{\mathcal{X}}_u(t)=(\mathcal{X}^i_u(t), i=1,\ldots,d)$ and $\vk{m}_u(t)=(m_u^i(t),i=1,\ldots,d)$.\\
Notice also that $g(\vk{v})=1$, and thus
\BQNY
|\vk{v}|=|\vk{v}|g(\vk{e})(g(\vk{e}))^{-1}=g(\vk{v})
(g(\vk{e}))^{-1}=(g(\vk{e}))^{-1}.
\EQNY
Now let us consider the probability
\begin{align*}
P_{u_k,\vk{e}}(y):=&\pk{\sup_{t\in[-S_1,S_2]}g(\vk{\mathcal{X}}_u(t)
+\vk{m}_u(t))>u_k}
=\pk{\sup_{t\in[-S_1,S_2]}g(\vk{e})|\vk{\mathcal{X}}_u(t)
+\vk{m}_u(t)|>u_k}\\
=&\pk{\sup_{t\in[-S_1,S_2]}\LT(|\vk{\mathcal{X}}_u(t)+\vk{m}_u(t)|^2
-(g(\vk{e}))^{-2}u_k^{2}\RT)>0}\\
=&\pk{\sup_{t\in[-S_1,S_2]}\LT(|\vk{\mathcal{X}}_u(t)+\vk{m}_u(t)|^2
-|\vk{v}|^2u_k^{2}\RT)>0}\\
=&\pk{\sup_{t\in[-S_1,S_2]}\LT(|\vk{\mathcal{X}}_u(t)+\vk{m}_u(t)
-\vk{v}u_k|^2
+2\LT\langle \vk{\mathcal{X}}_u(t)+\vk{m}_u(t)-\vk{v}u_k,\vk{v}u_k\RT\rangle\RT)>0}\\
=&\pk{\sup_{t\in[-S_1,S_2]}\zeta_u(t)>0},
\end{align*}
where $\zeta_u(t):=\LT\langle \vk{\mathcal{X}}_u(t)+\vk{m}_u(t)-\vk{v}u_k,\vk{v}u_k\RT\rangle
+\frac{1}{2}|\vk{\mathcal{X}}_u(t)+\vk{m}_u(t)-\vk{v}u_k|^2$.\\
As in \cite{PiterChaos2015}  we obtain
\BQNY
P_{u_k,\vk{e}}(y)\leq \mathbb{Q}_1\exp(-\mathbb{Q}_2y^2)+\mathbb{Q}_3\exp(-\mathbb{Q}_4y),
\EQNY
where $\mathbb{Q}_4>1$, and further
\BQN\label{bound}
P_{u_k}(y)\leq \mathbb{Q}_5\exp(-\mathbb{Q}_2y^2)+\mathbb{Q}_6\exp(-\mathbb{Q}_4y).
\EQN
By \eqref{uk}, \eqref{mean} and \eqref{var}, as $u\to \IF$
\BQNY
\E{\zeta_u(t)}
\rw -\LT(f(t)+a|t|^\alpha\RT)|\vk{v}|^2-y|\vk{v}|^2,
\EQNY
holds uniformly for all $t\in[-S_1,S_2]$ and $k\in K_u$, and moreover
\BQNY
\Var\LT(\zeta_u(t)-\zeta_u(s)\RT)\rw 2a|t-s|^\alpha|\vk{v}|^2,
\EQNY
holds uniformly for all $s,t\in[-S_1,S_2]$ and $k\in K_u$.
Given $B_{\alpha}(t),\ B^i_{\alpha}(t),\ i=1,2,\ldots,d$  independent fractional Brownian motion with the same Hurst index $\alpha/2\in(0,1]$ we can write
\begin{align*}
&\pk{\sup_{t\in[-S_1,S_2]}\sum_{i=1}^d\sqrt{2a}B_{\alpha}^i(t)v_i
-\LT(f(t)+a|t|^\alpha\RT)|\vk{v}|^2-y|\vk{v}|^2>0}\\
=&\pk{\sup_{t\in[-S_1,S_2]}\sqrt{2a}B_{\alpha}(t)|\vk{v}|
-\LT(f(t)+a|t|^\alpha\RT)|\vk{v}|^2>y|\vk{v}|^2}\\
=&\pk{\sup_{t\in[-S_1,S_2]}\sqrt{\frac{2a}{(g(\vk{e}))^{2}}}B_{\alpha}(t)
-\frac{1}{(g(\vk{e}))^{2}}f(t)-\frac{a}{(g(\vk{e}))^{2}}|t|^\alpha>y(g(\vk{e}))^{-2}}
=:P_{\vk{e}}(y).
\end{align*}
Consequently,
\BQNY
\lim_{u\rw\IF}\underset{y\in[-W,W]}{\sup_{k\in K_u}}\LT|P_{u_k,\vk{e}}(y)-P_{\vk{e}}(y)\RT|=0.
\EQNY
By \cite{PiterChaos2015}, we know that there exist a density function $j(\vk{e})$ on $\mathcal{M}$ such that
\BQNY
\lim_{u\rw\IF}\sup_{k\in K_u }\LT|\int_{\mathbb{S}_{d-1}}f_k(\vk{e})P_{u_k,\vk{e}}(y)d \vk{e}-\int_{\mathcal{M}}j(\vk{e})P_{\vk{e}}(y)d \vk{e}\RT|=0,
\EQNY
and
\BQNY
\int_{\mathcal{M}}j(\vk{e})P_{\vk{e}}(y)d \vk{e}
=\int_{\mathcal{M}}j(\vk{e})\pk{\sup_{t\in[-S_1,S_2]}\eta(t)>y} d\vk{e}=\pk{\sup_{t\in[-S_1,S_2]}\eta(t)>y},
\EQNY
where we used the fact that $g(\vk{e})=1$ for $\vk{e}\in \mathcal{M}$ and $\int_{\mathcal{M}}j(\vk{e})d\vk{e}=1$. Hence we have
\BQN\label{pp}
\lim_{u\rw\IF}\sup_{k\in K_u, y\in[-W,W]}\LT|P_{u_k}(y)-\pk{\sup_{t\in[-S_1,S_2]}\eta(t)>y}\RT|=0.
\EQN
{\bf Step 2}:
For $\vn_1<\frac{1}{2}$, we have (below we set $U_u(t):= g(\vk{\xi}(u^{-2/\alpha}t+t_0))$)
\begin{align*}
I_2(u)\leq&u_k^{-1}\int_{W}^{\IF}\pk{\sup_{t\in[-S_1,S_2]}
\ehe{U_u(t)})>u_k\Big|g(\vk{\xi}(t_0))=\uy}
\ehe{p}(\uy)dy\\
=&u_k^{-1}\int_{W}^{u_k^{2-\vn_1}}\pk{\sup_{t\in[-S_1,S_2]}
\ehe{U_u(t)})>u_k\Big|g(\vk{\xi}(t_0))=\uy}
\ehe{p}(\uy)dy\\
&+u_k^{-1}\int_{u_k^{2-\vn_1}}^{\IF}\pk{\sup_{t\in[-S_1,S_2]}
\ehe{U_u(t)})>u_k
\Big|g(\vk{\xi}(t_0))=\uy}
\ehe{p}(\uy)dy\\
=&u_k^{-1}\int_{W}^{u_k^{2-\vn_1}}\pk{\sup_{t\in[-S_1,S_2]}
\ehe{U_u(t)})>u_k\Big|g(\vk{\xi}(t_0))=\uy}
\ehe{p}(\uy)dy\\
&+
\pk{\sup_{t\in[-S_1,S_2]} U_u(t)
>u_k,\ g(\vk{\xi}(t_0))\leq u_k- u_k^{1-\vn_1}}\\
=&: J_1(u)+J_2(u).
\end{align*}
Then we use the similar argumentation in {\bf Step 1} to deal with
$$P_{u_k}(y)=\pk{\sup_{t\in[-S_1,S_2]}\ehe{U_u(t)})
>u_k\Big|g(\vk{\xi}(t_0))=\uy},\ y\in [W,u_k^{2-\vn_1}],$$
and we get
\BQNY
\lim_{u\rw\IF}\sup_{k\in K_u, y\in[W,u_k^{2-\vn_1}]}\LT|P_{u_k}(y)-\pk{\sup_{t\in[-S_1,S_2]}\eta_1(t)>y}\RT|=0,
\EQNY
where $\eta_1(t)=\sqrt{2a}B_\alpha(t)-a\abs{t}^\alpha$.
\COM{Similar with \eqref{bb}, we drive that
\BQNY
&&\sup_{k\in K_u, y\in[W,u_k^{\frac{2-\vn_1}{p}}]}\LT|
\frac{u_k^{1-2/p}\ehe{p}(u_k-u_k^{1-2/p}y)}
{\frac{1}{p}\pk{g(\vk{\xi}(t_0))>u_k}e^{\frac{y}{p}}}-1\RT|\nonumber\\
&&=\sup_{k\in K_u, y\in[W,u_k^{\frac{2-\vn_1}{p}}]}\LT|\frac{u_k^{1-2/p}\frac{h_0g^{-(m+1)/p}}{p}
(1+R(u_k-u_k^{1-2/p}y))\LT(u_k-u_k^{1-2/p}y\RT)^{\frac{m+1}{p}-1}
\exp\LT(-\frac{1}{2}\LT(\frac{u_k-u_k^{1-2/p}y}{g}\RT)^{2/p}\RT)}{\frac{1}{p}
\pk{g(\vk{\xi}(t_0))>u_k}e^{\frac{y}{p}}}-1\RT|\nonumber\\
&&=\sup_{k\in K_u, y\in[W,u_k^{\frac{2-\vn_1}{p}}]}\LT|\frac{(1+R(u_k-u_k^{1-2/p}y))\LT(1-u_k^{-2/p}y\RT)
^{\frac{m+1}{p}-1}\exp\LT(-\frac{\theta(u_k,y)}{2u_k^{2/p}}\RT)}
{1+R_1(u_k)}-1\RT|\nonumber\\
&&\rw 0,\  u\rw\IF.
\EQNY}
Then by \eqref{bb} and similarly to \eqref{boundeta}
\begin{align}\label{ll}
\sup_{k\in K_u}\frac{J_1(u)}{\pk{g(\vk{\xi}(t_0))>u_k}}&\leq \mathbb{Q}_7 \int_{W}^{\IF}\pk{\sup_{t\in[-S_1,S_2]}\eta_1(t)>y}e^{y}dy\nonumber\\
&\leq \mathbb{Q}_8\int_{W}^{\IF}e^{y} y^{2/\alpha}\Psi\LT(\frac{y}{\sqrt{2a\max(S_1,S_2)}}\RT) d y
\rw 0, \ W\rw\IF.
\end{align}
Now let us proceed the analysis of $J_2(u)$. We have
\BQNY
J_2(u)\leq\pk{\sup_{t\in[-u^{-2/\alpha}S_1,u^{-2/\alpha}S_2]}\LT(g(\vk{\xi}(t+t_0))
-g(\vk{\xi}(t_0))\RT)> u_k^{1-\vn_1}}.
\EQNY
Since
\BQNY
g(\vk{x}+\vk{y})-g(\vk{x})=\int_{0}^1\langle\nabla g(\vk{x}+h\vk{y}),\vk{y}\rangle d h,
\EQNY
and
\BQNY
\LT|g(\vk{x}+\vk{y})-g(\vk{x})\RT|\leq g_1|\vk{y}|,
\EQNY
where $g_1=\max_{|\vk{v}|=1}\LT|\nabla g(\vk{v})\RT|$. Hence, denoting $\Delta\vk{\xi}(t)=\vk{\xi}(t+t_0)-\vk{\xi}(t_0)$, we get that
\BQNY
\sup_{t\in[-u^{-2/\alpha}S_1,u^{-2/\alpha}S_2]}
\LT(g(\vk{\xi}(t+t_0))-g(\vk{\xi}(t_0))\RT)\leq g_1 \sup_{t\in[-u^{-2/\alpha}S_1,u^{-2/\alpha}S_2]}\LT|\Delta\vk{\xi}(t)\RT|
.
\EQNY
Further, since by \eqref{uk} for $\vn>0$ when $u$ large enough
\BQN\label{uu}
\inf_{k\in K_u}u_k\geq(1-\vn)u,
\EQN
then we have
\BQNY
&&\pk{\sup_{t\in[-S_1,S_2]}\LT(\ehe{U_u(t)})-g(\vk{\xi}(t_0))\RT)> u_k^{1-\vn_1}}
\leq\pk{\sup_{t\in[-u^{-2/\alpha}S_1,u^{-2/\alpha}S_2]}\LT(g(\vk{\xi}(t+t_0))
-g(\vk{\xi}(t_0))\RT)> u_k^{1-\vn_1}}\\
&&\leq\pk{\sup_{t\in[-u^{-2/\alpha}S_1,u^{-2/\alpha}S_2]}
\LT|t^{-\alpha/2}\Delta\vk{\xi}(t)\RT|
> \mathbb{Q}_7u^{2-\vn_1}},
\EQNY
where $\mathbb{Q}_7=(1-\vn)^{1-\vn_1}g_1^{-1}(\max(S_1,S_2))^{-\alpha/2}$.
Now consider the vector process
\BQNY
\vk{X}_1(t):=t^{-\alpha/2}\Delta\vk{\xi}(t),
t\in[-u^{-2/\alpha}S_1,u^{-2/\alpha}S_2],
\EQNY
which tends weakly to zero, hence
\BQNY
\pk{\sup_{t\in[-S_1,S_2]}\LT(\ehe{U_u(t)}-g(\vk{\xi}(t_0))\RT)> u_k^{1-\vn_1}}
\leq\pk{\sup_{t\in[-u^{-2/\alpha}S_1,u^{-2/\alpha}S_2]}G(\vk{X}_1(t))
>\mathbb{Q}_7u^{2-\vn_1}},
\EQNY
where $G(\vk{x})=|\vk{x}|,\ \vk{x}\in\R^{d}$, is a homogeneous function of order $1$. Now from \netheo{Cor2} it follows that the probability $J_2(u)$ is exponentially smaller than the probability $I_1(u)$.

 ii) Next, we prove \eqref{P2p2}.  By \eqref{lrr11}, for $\vn_T\in (0,a)$ we can find $T$ small enough such that
\BQN \label{abound}
(a-\vn_T)\leq a(t)\leq (a+\vn_T),\quad (a-\vn_T) \abs{t-s}^\alpha \leq 1-r(s,t)\leq (a+\vn_T) \abs{t-s}^\alpha
\EQN
holds for $s,t\in [t_0-T,t_0+T]$.
Set $$\vk{\xi}'_{(u,k)}(t)=
\vk{\xi}(u^{-2/\alpha}(a(ku^{-2/\alpha}S+t_0))^{-1/\alpha}t+ku^{-2/\alpha}S+t_0),
\quad t\in[0,(a(ku^{-2/\alpha}S+t_0))^{1/\alpha}S].$$
Then similar to i), we have for any $W>0$ and all $u$ large
\begin{align*}
&\pk{\sup_{t\in[0,S]}g(\vk{\xi}(u^{-2/\alpha}t+ku^{-2/\alpha}S+t_0))>u_k}\\
&=\int_{-\IF}^{\IF}\pk{\sup_{t\in[0,(a(ku^{-2/\alpha}S+t_0))^{1/\alpha}S]}
g(\vk{\xi}'_{(u,k)}(t))>u_k\Big|g(\vk{\xi}'_{(u,k)}(0))=x}
\ehe{p}(x)dx\\
&=u_k^{-1}\int_{-\IF}^{\IF}\pk{\sup_{t\in[0,(a(ku^{-2/\alpha}S+t_0))^{1/\alpha}S]}
g(\vk{\xi}'_{(u,k)}(t))>u_k\Big|g(\vk{\xi}'_{(u,k)}(0))=\uy}
\ehe{p}(\uy)dy\\
&=u_k^{-1}\int_{-W}^{W}\pk{\sup_{t\in[0,(a(ku^{-2/\alpha}S+t_0))^{1/\alpha}S]}
g(\vk{\xi}'_{(u,k)}(t))>u_k\Big|g(\vk{\xi}'_{(u,k)}(0))=\uy}
\ehe{p}(\uy)dy\\
&\quad+\pk{\sup_{t\in[0,(a(ku^{-2/\alpha}S+t_0))^{1/\alpha}S]}
g(\vk{\xi}'_{(u,k)}(t))>u_k,\ g(\vk{\xi}'_{(u,k)}(0))\leq u_k-u_k^{-1}W}\\
&\quad+\pk{\sup_{t\in[0,(a(ku^{-2/\alpha}S+t_0))^{1/\alpha}S]}
g(\vk{\xi}'_{(u,k)}(t))>u_k, g(\vk{\xi}'_{(u,k)}(0))> u_k+u_k^{-1}W}\\
&=: I'_1(u)+I'_2(u)+I'_3(u),
\end{align*}
where we used the fact that
\BQNY
\ehe{p}(x)=p_{g(\vk{\xi}'_{(u,k)}(0))}(x),\quad x\in\R.
\EQNY
By
\BQNY
I'_3(u)\leq \pk{g(\vk{\xi}'_{(u,k)}(0))> u_k+u_k^{-1}W}=
\pk{g(\vk{\xi}(t_0))> u_k+u_k^{-1}W},
\EQNY
we know that as in \eqref{II3}
\BQNY
\sup_{k\in K_u}\frac{I'_3(u)}{\pk{g(\vk{\xi}(t_0))>u_k}}\leq \sup_{k\in K_u} \frac{\pk{g(\vk{\xi}(t_0))> u_k+u_k^{-1}W}}{\pk{g(\vk{\xi}(t_0))>u_k}} \rw 0, \ u\rw\IF, W\rw \IF.
\EQNY
Next we consider $I'_1(u)$. Setting
\BQNY
P'_{u_k}(y)=\pk{\sup_{t\in[0,(a(ku^{-2/\alpha}S+t_0))^{1/\alpha}S]}
g(\vk{\xi}'_{(u,k)}(t))>u_k\Big|g(\vk{\xi}'_{(u,k)}(0))=\uy},
\EQNY
we have by \eqref{abound}
\BQNY
P'_{u_k}(y)\geq \pk{\sup_{t\in[0,(a-\vn_T)^{1/\alpha}S]}
g(\vk{\xi}'_{(u,k)}(t))>u_k\Big|g(\vk{\xi}'_{(u,k)}(0))=\uy},
\EQNY
and
\BQNY
P'_{u_k}(y)
\leq \pk{\sup_{t\in[0,(a+\vn_T)^{1/\alpha}S]}
g(\vk{\xi}'_{(u,k)}(t))>u_k\Big|g(\vk{\xi}'_{(u,k)}(0))=\uy}.
\EQNY
Similar to the proof of i), we analyse
\BQNY
I^{+}_u(W)=u_k^{-1}\int_{-W}^{W}\pk{\sup_{t\in[0,(a+\vn_T)^{1/\alpha}S]}
g(\vk{\xi}'_{(u,k)}(t))>u_k\Big|g(\vk{\xi}'_{(u,k)}(0))=\uy}
\ehe{p}(\uy)dy
\EQNY
and
\BQNY
I^{-}_u(W)=u_k^{-1}\int_{-W}^{W}\pk{\sup_{t\in[0,(a-\vn_T)^{1/\alpha}S]}
g(\vk{\xi}'_{(u,k)}(t))>u_k\Big|g(\vk{\xi}'_{(u,k)}(0))=\uy}
\ehe{p}(\uy)dy
\EQNY
for some constant $W>0$.
Next we use the similar arguments as in ii) to analyse $I^{+}_u(W)$ with
$\eta_2(t)=\sqrt{2}B_{\alpha}(t)-\abs{t}^\alpha$.
Since for $t\in [0,(a+\vn_T)^{1/\alpha}S]$ and all $k\in K_u$
\BQNY
Cov(\xi'_{(u,k),i}(t),\xi'_{(u,k),i}(0))\sim 1-a(ku^{-2/\alpha}S+t_0)\abs{u^{-2/\alpha}
(a(ku^{-2/\alpha}S+t_0))^{-1/\alpha}t}^\alpha
=1-u^{-2}\abs{t}^\alpha,
\EQNY
as $u\rw\IF$, we replace \eqref{mean} and \eqref{var} with
\BQNY
m_u^i(t)=u_k v_i\LT(1-|t|^\alpha u^{-2}-y u_k^{-2}+u^{-2}\rho_1(u,k,t,y)\RT),\quad
V_u(s,t)=2u^{-2}(|t-s|^\alpha+\rho_2(u,s,t)).
\EQNY
Further for $\zeta_u(t):=\LT\langle \vk{\mathcal{X}}_u(t)+\vk{m}_u(t)-\vk{v}u_k,\vk{v}u_k\RT\rangle
+\frac{1}{2}|\vk{\mathcal{X}}_u(t)+\vk{m}_u(t)-\vk{v}u_k|^2$
where $\vk{\mathcal{X}}_u(t)$ is Gaussian vector process where $\mathcal{X}^i_u(t), i=1,2,\ldots,d$ have zero means, continuous trajectories with $t\in[0,(a+\vn_T)^{1/\alpha}S]$, and covariance functions equal to the respective conditional covariance function of $\xi'_{(u,k),i}(t)$. then
\BQNY
\E{\zeta_u(t)}
\rw -|t|^\alpha|\vk{v}|^2-y|\vk{v}|^2,
\EQNY
holds uniformly for all $t\in[0,(a+\vn_T)^{1/\alpha}S]$, $y\in[-W,W]$ and $k\in K_u$ and
\BQNY
\Var\LT(\zeta_u(t)-\zeta_u(s)\RT)\rw 2|t-s|^\alpha|\vk{v}|^2,
\EQNY
holds uniformly for all $s,t\in[0,(a+\vn_T)^{1/\alpha}S]$, $y\in[-W,W]$ and $k\in K_u$.
Then
\begin{align*}
&\pk{\sup_{t\in[0,(a+\vn_T)^{1/\alpha}S]}\sum_{i=1}^d\sqrt{2}B_{\alpha}^i(t)v_i
-|t|^\alpha|\vk{v}|^2-y|\vk{v}|^2>0}
=\pk{\sup_{t\in[0,(a+\vn_T)^{1/\alpha}S]}\sqrt{2}B_{\alpha}(t)|\vk{v}|
-|t|^\alpha|\vk{v}|^2>y|\vk{v}|^2}\\
&=\pk{\sup_{t\in[0,(a+\vn_T)^{1/\alpha}S]}\sqrt{\frac{2}{(g(\vk{e}))^{2}}}B_{\alpha}(t)
-\frac{1}{(g(\vk{e}))^{2}}|t|^\alpha>y(g(\vk{e}))^{-2}}
=:P^{+}_{\vk{e}}(y),
\end{align*}
and similarly we analyse  $I^{-}_u(W)$ with
\BQNY
P^{-}_{\vk{e}}(y):=\pk{\sup_{t\in[0,(a-\vn_T)^{1/\alpha}S]}\sqrt{\frac{2}
{(g(\vk{e}))^{2}}}B_{\alpha}(t)
-\frac{1}{(g(\vk{e}))^{2}}|t|^\alpha>y(g(\vk{e}))^{-2}}.
\EQNY
Thus we have that for $u$ large enough
\BQNY
\pk{\sup_{t\in[0,(a-\vn_T)^{1/\alpha}S]}\eta_2(t)>y}\leq P'_{u_k}(y)
\leq\pk{\sup_{t\in[0,(a+\vn_T)^{1/\alpha}S]}\eta_2(t)>y}
\EQNY
 holds for all $k\in K_u$ and $y\in[-W,W]$, which combining with \eqref{bb} drives that for $u$ large enough and any $W>0$
 \BQNY
\mathcal{H}_{\alpha}[0,(a-\vn_T)^{1/\alpha}S]\leq \frac{I'_1(u)}{\pk{g(\vk{\xi}(t_0))>u_k}}\leq\mathcal{H}_{\alpha}[0,(a+\vn_T)^{1/\alpha}S],
 \EQNY
 holds for any $k\in K_u $.\\
Now we consider $I'_2(u)$. Similar to {\bf Step 2} of i), for $\vn_1<\frac{1}{2}$, we have
\begin{align*}
I_2(u)&\leq u_k^{-1}\int_{W}^{\IF}\pk{\sup_{t\in[0,(a(ku^{-2/\alpha}S+t_0))^{1/\alpha}S]}
g(\vk{\xi}'_{(u,k)}(t))>u_k\Big|g(\vk{\xi}'_{(u,k)}(0))=\uy}
\ehe{p}(\uy)dy\\
&=u_k^{-1}\int_{W}^{u_k^{2-\vn_1}}
\pk{\sup_{t\in[0,(a(ku^{-2/\alpha}S+t_0))^{1/\alpha}S]}
g(\vk{\xi}'_{(u,k)}(t))>u_k\Big|g(\vk{\xi}'_{(u,k)}(0))=\uy}
\ehe{p}(\uy)dy\\
&\quad+u_k^{-1}\int_{u_k^{2-\vn_1}}^{\IF}
\pk{\sup_{t\in[0,(a(ku^{-2/\alpha}S+t_0))^{1/\alpha}S]}
g(\vk{\xi}'_{(u,k)}(t))>u_k\Big|g(\vk{\xi}'_{(u,k)}(0))=\uy}
\ehe{p}(\uy)dy\\
&=u_k^{-1}\int_{W}^{u_k^{2-\vn_1}}
\pk{\sup_{t\in[0,(a(ku^{-2/\alpha}S+t_0))^{1/\alpha}S]}
g(\vk{\xi}'_{(u,k)}(t))>u_k\Big|g(\vk{\xi}'_{(u,k)}(0))=\uy}
\ehe{p}(\uy)dy\\
&\quad+
\pk{\sup_{t\in[0,(a(ku^{-2/\alpha}S+t_0))^{1/\alpha}S]}g(\vk{\xi}'_{(u,k)}(t))>u_k,\ g(\vk{\xi}'_{(u,k)}(0))\leq u_k- u_k^{1-\vn_1}}\\
&=: J'_1(u)+J'_2(u).
\end{align*}
Then we use the similar argumentation in the former step to deal with
$$P'_{u_k}(y)=\pk{\sup_{t\in[0,(a(ku^{-2/\alpha}S+t_0))^{1/\alpha}S]}
g(\vk{\xi}'_{(u,k)}(t))>u_k\Big|g(\vk{\xi}'_{(u,k)}(0))=\uy},
 \ y\in [W,u_k^{2-\vn_1}],$$
and we get
\BQNY
\lim_{u\rw\IF}\sup_{k\in K_u, y\in[W,u_k^{2-\vn_1}]}P'_{u_k}(y)
\leq\pk{\sup_{t\in[0,(a+\vn_T)^{1/\alpha}S]}\eta_2(t)>y}.
\EQNY
With \eqref{bb} and similarly to \eqref{boundeta}, we derive that
\begin{align*}
\sup_{k\in K_u}\frac{J'_1(u)}{\pk{g(\vk{\xi}(t_0))>u_k}}&\leq \mathbb{Q}_{11} \int_{W}^{\IF}\pk{\sup_{t\in[0,(a+\vn_T)^{1/\alpha}S]}
\eta_2(t)>y}e^{y}dy\nonumber\\
&\leq \mathbb{Q}_{12} \int_{W}^{\IF}e^{y}y^{2/\alpha}\Psi\LT(\frac{y}{\sqrt{2(a+\vn_T)^{1/\alpha}S}}\RT) d y \rw 0, \ W\rw\IF.
\end{align*}
Since
\begin{align*}
J'_2(u)&=\pk{\sup_{t\in[0,(a(ku^{-2/\alpha}S+t_0))^{1/\alpha}S]}g(\vk{\xi}'_{(u,k)}(t))>u_k,\ g(\vk{\xi}'_{(u,k)}(0))\leq u_k- u_k^{1-\vn_1}}\\
&\leq\pk{\sup_{t\in[0,S]}g(\vk{\xi}(u^{-2/\alpha}t
+ku^{-2/\alpha}S+t_0))>u_k,\ g(\vk{\xi}(
ku^{-2/\alpha}S+t_0))\leq u_k- u_k^{1-\vn_1}},
\end{align*}
then using the same argument for $J_2(u)$ as in {\bf Step 2} of i) with $t_0$ replaced by $ku^{-2/\alpha}S+t_0$, we obtain that
\BQNY
\lim_{u\rw\IF}\sup_{k\in K_u}\frac{J'_2(u)}{\pk{g(\vk{\xi}(t_0))>u_k}}=0.
\EQNY
Thus we finish the proof  of \eqref{P2p2}.\\
If we let $T\rw 0$ in \eqref{P2p2}, then $\vn_T\rw 0$ and \eqref{P2p} follows.

\QED

\prooflem{in1}
We use several results and arguments from ii) of the proof of \nelem{im1}.
Set
\BQNY
&&\vk{\xi}'_{(u,k)}(t)=
\vk{\xi}(u^{-2/\alpha}(a(ku^{-2/\alpha}S+t_0))^{-1/\alpha}t+ku^{-2/\alpha}S+t_0),
\quad t\in\R,\quad a_{u,k}=(a(ku^{-2/\alpha}S+t_0))^{1/\alpha},\\
&&\mathcal{A}'_i(u_k)=\LT\{\sup_{t\in [a_{u,k}T_i,a_{u,k}(T_i+S)]}
g(\vk{\xi}'_{(u,k)}(t))>u_k\RT\},\quad i=1,2.
\EQNY
For any $W>0$ and all $u$ large
\begin{align*}
\pk{\mathcal{A}_1(u_k),\mathcal{A}_2(u_k)}
&=\int_{-\IF}^{\IF}\pk{\mathcal{A}'_1(u_k),\mathcal{A}'_2(u_k)
\Big|g(\vk{\xi}'_{(u,k)}(0))=x}
\ehe{p}(x)dx\\
&=u_k^{-1}\int_{-\IF}^{\IF}\pk{\mathcal{A}'_1(u_k),\mathcal{A}'_2(u_k)
\Big|g(\vk{\xi}'_{(u,k)}(0))=\uy}
\ehe{p}(\uy)dy\\
&=u_k^{-1}\int_{-W}^{W}\pk{\mathcal{A}'_1(u_k),\mathcal{A}'_2(u_k)
\Big|g(\vk{\xi}'_{(u,k)}(0))=\uy}
\ehe{p}(\uy)dy\\
&\quad+\pk{\mathcal{A}'_1(u_k),\mathcal{A}'_2(u_k),\ g(\vk{\xi}'_{(u,k)}(0))\leq u_k-u_k^{-1}W}\\
&\quad+\pk{\mathcal{A}'_1(u_k),\mathcal{A}'_2(u_k),g(\vk{\xi}'_{(u,k)}(0))> u_k+u_k^{-1}W}\\
&=: I_1(u)+I_2(u)+I_3(u),
\end{align*}
where $\uy=u_k-u_k^{-1}y$ and  $\ehe{p}(x)$ is the density function of $g(\vk{\xi}(t_0))$ which is showed in \netheo{Them00}.
Similar to the proof of \nelem{im1}  ii), we know
\BQNY
\sup_{k\in K_u}\frac{I_2(u)}{\pk{g(\vk{\xi}(t_0))>u_k}}\leq \sup_{k\in K_u}\frac{\pk{\mathcal{A}'_1(u_k),\ g(\vk{\xi}'_{(u,k)}(0))\leq u_k-u_k^{-1}W}}{\pk{g(\vk{\xi}(t_0))>u_k}}\rw 0, \ u\rw\IF.
\EQNY
Further, as in \eqref{II3}
\BQNY
\sup_{k\in K_u}\frac{I_3(u)}{\pk{g(\vk{\xi}(t_0))>u_k}}\leq \sup_{k\in K_u}\frac{\pk{ g(\vk{\xi}(t_0))\leq u_k-u_k^{-1}W}}{\pk{g(\vk{\xi}(t_0))>u_k}}\rw 0,\ u\rw\IF,\ W\rw \IF,
\EQNY
then we can choose $W$ large enough such that $\frac{I_3(u)}{\pk{g(\vk{\xi}(t_0))>u_k}}$ less that $\mathbb{Q}_1\exp\LT(-\frac{a}{8}|T_2-T_1-S|^\alpha\RT)$.\\
Set below
\BQNY
&&\eta(t)=\sqrt{2}B_{\alpha}(t)-\abs{t}^\alpha,\quad \widetilde{P}_{u_k}(y,z)=\pk{\widetilde{\mathcal{A}}'_1(u_k,z),\widetilde{\mathcal{A}}'_2(u_k,z)
\Big|g(\vk{\xi}'_{(u,k)}(0))=\uy},\\
&&\widetilde{\mathcal{A}}'_i(u_k,z)=\LT\{\sup_{t\in [zT_i,z(T_i+S)]}
g(\vk{\xi}'_{(u,k)}(t))>u_k\RT\},\quad z\in[\underline{a},\overline{a}], \quad i=1,2,
\EQNY
and  $\underline{a}=\inf_{t\in[t_0-\vn_0,t_0+\vn_0]}a(t)$ and $\overline{a}=\sup_{t\in[t_0-\vn_0,t_0+\vn_0]}a(t)$.
First note that
\BQN\label{upppp1}
I_1(u)\leq \sup_{z\in[\underline{a},\overline{a}]}u_k^{-1}\int_{-W}^{W}\widetilde{P}_{u_k}(y,z)
\ehe{p}(\uy)dy.
\EQN
Similarly to the arguments as \eqref{basic} in proof of \nelem{im1}
we have that for $u$ large enough
\begin{align*}
&\underset{z\in[\underline{a},\overline{a}]}{\sup_{k\in K_u}}\LT|\frac{u_k^{-1}\int_{-W}^{W}\widetilde{P}_{u_k}(y,z)
\ehe{p}(\uy)dy}
{\pk{g(\vk{\xi}(t_0))>u_k}}-\int_{-W}^{W}e^{y}\pk{\bigcap_{i=1,2}\LT(\sup_{t\in[zT_i,
z(T_i+S)]}\eta(t)>y\RT)} d y\RT|\nonumber\\
&\leq\underset{z\in[\underline{a},\overline{a}]}{\sup_{k\in K_u}} e^{W}\int_{-W}^{W}\LT|
\frac{u_k^{-1}\ehe{p}(\uy)}
{e^{y}\pk{g(\vk{\xi}(t_0))>u_k}}-1\RT|d y \\
&\quad +\underset{z\in[\underline{a},\overline{a}]}{\sup_{k\in K_u}}e^{W}\int_{-W}^{W}\LT|\widetilde{P}_{u_k}(y,z)-\pk{\bigcap_{i=1,2}\LT(\sup_{t\in[zT_i,
z(T_i+S)]}\eta(t)>y\RT)}\RT|d y,
\end{align*}
where the first term goes to zero as $u\rw\IF$ by \eqref{bb}.\\
By similar arguments as in {\bf Step 2} of the \nelem{im1} we have for any $(y,z)\in[-W,W]\times[\underline{a},\overline{a}]$
\BQNY
c_u(y,z):={\sup_{k\in K_u }}\abs{\widetilde{P}_{u_k}(y,z)-\pk{\bigcap_{i=1,2}\LT(\sup_{t\in[z T_i,
z(T_i+S)]}\eta(t)>y\RT)}}\rw 0,\quad u\rw\IF.
\EQNY
By the well-known  Severini-Egorov theorem, for any $\vn>0,W>0$ the convergence
$$ c_u(y,z) \to 0, \quad u\to \IF$$
is uniform  for $(y,z)\in\LT([-W,W]\times[\underline{a},\overline{a}]\RT)\setminus K_\vn$ with $K_\vn$ a measurable set with Lebesgue measure not exceeding $\vn$. Hence, we can write
\BQNY
\lim_{u\rw\IF}\underset{(y,z)\in\LT([-W,W]\times[\underline{a},\overline{a}]\RT)\setminus K_\vn}{\sup_{k\in K_u }}\abs{\widetilde{P}_{u_k}(y,z)-\pk{\bigcap_{i=1,2}\LT(\sup_{t\in[z T_i,
z(T_i+S)]}\eta(t)>y\RT)}}=0.
\EQNY

\COM{Then we have
\BQNY
\lim_{u\rw\IF}{\sup_{k\in K_u }}\abs{\frac{I_1(u)}{\pk{g(\vk{\xi}(t_0))>u_k}}- \int_{-{W}/{p\gg^{2/p}}}^{{W}/{p\gg^{2/p}}}e^{v}
\pk{\bigcap_{i=1,2}\LT(\sup_{t\in[a_{u,k}T_i,a_{u,k}(T_i+S)]}\eta(t)>x\RT)}dv}=0.
\EQNY}
Since $\vn>0$ can be chosen arbitrary small, we obtain
\BQNY
\lim_{u\rw\IF}\underset{z\in[\underline{a},\overline{a}]}{\sup_{k\in K_u}}e^{W}\int_{-W}^{W}\LT|\widetilde{P}_{u_k}(y,z)-\pk{\bigcap_{i=1,2}\LT(\sup_{t\in[zT_i,
z(T_i+S)]}\eta(t)>y\RT)}\RT|d y=0.
\EQNY
Thus we have as $u\rw\IF$
\BQN\label{ASY1}
\underset{z\in[\underline{a},\overline{a}]}{\sup_{k\in K_u}}\LT|\frac{u_k^{-1}\int_{-W}^{W}\widetilde{P}_{u_k}(y,z)
\ehe{p}(\uy)dy}
{\pk{g(\vk{\xi}(t_0))>u_k}}-\int_{-W}^{W}e^{y}\pk{\bigcap_{i=1,2}\LT(\sup_{t\in[zT_i,
z(T_i+S)]}\eta(t)>y\RT)} d y\RT|\rw0.
\EQN
By Slepian inequality in \cite{Pit96} and \cite{Uniform2016}[Theorem 3.1], for all $u$ large
\BQNY
\underset{z\in[\underline{a},\overline{a}]}{\sup_{k\in K_u}}\frac{1}{\Psi(u_k)}\pk{\bigcap_{i=1,2}\LT(\sup_{t\in[zT_i,
z(T_i+S)]}\xi'_{(u,k),1}(t)>u_k\RT)}
\leq\mathbb{Q}_2\exp\LT(-\frac{\underline{a}}{8}\abs{T_2-T_1-S}^\alpha\RT).
\EQNY

Further, we have for $W>0$
\BQNY
&&\pk{\bigcap_{i=1,2}\LT(\sup_{t\in[zT_i,
z(T_i+S)]}\xi'_{(u,k),1}(t)>u_k\RT)}\\
&&=\frac{1}{\sqrt{2\pi u_k}}\int_{-\IF}^{\IF}e^{-\frac{1}{2}\LT(\uy\RT)^2}
\pk{\bigcap_{i=1,2}\LT(\sup_{t\in[zT_i,
z(T_i+S)]}\xi'_{(u,k),1}(t)>u_k\RT)
\Big|\xi'_{(u,k),1}(0)=\uy}dy\\
&&=\frac{1}{\sqrt{2\pi u_k}}\int_{-\IF}^{\IF}e^{-\frac{1}{2}\LT(\uy\RT)^2}
\pk{\bigcap_{i=1,2}\LT(\sup_{t\in[zT_i,
z(T_i+S)]}\mathcal{X}_{u,k}^y(t)>y\RT)}dy\\
&&\geq \Psi(u_k)\int_{-W}^{W}e^{y}
\pk{\bigcap_{i=1,2}\LT(\sup_{t\in[zT_i,
z(T_i+S)]}\mathcal{X}_{u,k}^y(t)>y\RT)}dy,
\EQNY
where
$$\mathcal{X}_{u,k}^y(t)=u_k\LT(\xi'_{(u,k),1}(t)-u_k\RT)+y
\Big|\xi'_{(u,k),1}(0)=\uy.$$
Thus we have
\BQN\label{RR1}
\underset{z\in[\underline{a},\overline{a}]}{\sup_{k\in K_u}}\int_{-W}^{W}e^{y}
\pk{\bigcap_{i=1,2}\LT(\sup_{t\in[zT_i,
z(T_i+S)]}\mathcal{X}_{u,k}^y(t)>y\RT)}dy\leq
\mathbb{Q}_2\exp\LT(-\frac{\underline{a}}{8}\abs{T_2-T_1-S}^\alpha\RT).
\EQN
From condition \eqref{lrr11}, we obtain uniformly with respect to $t\in [zT_1,
z(T_1+S)] \cup[zT_2,
z(T_2+S)], z\in[\underline{a},\overline{a}]$, $y\in[-W,W]$ and $k\in K_u$ that
\BQN
\E{\mathcal{X}^{y}_{u,k}(t)}\rw -|t|^\alpha,\quad u\rw\IF, \label{XM12}
\EQN
and also for any $t,t'\in [zT_1,
z(T_1+S)] \cup[zT_2,
z(T_2+S)],z\in[\underline{a},\overline{a}]$, $y\in\R$ and $k\in K_u$
\BQN\label{XVar}
\Var\LT(\mathcal{X}^{y}_{u,k}(t)-\mathcal{X}^{y}_{u,k}(t')\RT)\rw 2|t-t'|^\alpha,\quad u\rw\IF.
\EQN
 Consequently, since supremum is a continuous functional, for any $y\in[-W,W]$, we have
\BQNY
\lim_{u\rw\IF}\underset{y\in[-W,W],z\in[\underline{a},\overline{a}]}{\sup_{k\in K_u}} \abs{\pk{\bigcap_{i=1,2}\LT(\sup_{t\in[zT_i,
z(T_i+S)]}\mathcal{X}_{u,k}^y(t)>y\RT)}- \pk{\bigcap_{i=1,2}\LT(\sup_{t\in[zT_i,
z(T_i+S)]}\eta(t)>y\RT)}}=0.
\EQNY
The above equality implies
\BQNY
\lim_{u\rw\IF}\underset{z\in[\underline{a},\overline{a}]}{\sup_{k\in K_u}} \abs{\int_{-W}^{W}e^y\pk{\bigcap_{i=1,2}\LT(\sup_{t\in[zT_i,
z(T_i+S)]}\mathcal{X}_{u,k}^y(t)>y\RT)}dy- \int_{-W}^{W}e^y\pk{\bigcap_{i=1,2}\LT(\sup_{t\in[zT_i,
z(T_i+S)]}\eta(t)>y\RT)}dy}=0,
\EQNY
which combined with \eqref{upppp1}, \eqref{ASY1} and \eqref{RR1} yields
\BQN
\frac{I_1(u)}{\pk{g(\vk{\xi}(t_0))>u_k}}\leq \mathbb{Q}_2\exp\LT(-\frac{\underline{a}}{8}\abs{T_2-T_1-S}^\alpha\RT)
\EQN
establishing the proof.
\COM{and
\BQNY
\lim_{u\rw\IF}{\sup_{k\in K_u }}\abs{\frac{I_1(u)}{\pk{g(\vk{\xi}(t_0))>u_k}}- \int_{-{W}/{p\gg^{2/p}}}^{{W}/{p\gg^{2/p}}}e^{v}
\pk{\bigcap_{i=1,2}\LT(\sup_{t\in[a_{u,k}T_i,a_{u,k}(T_i+S)]}\mathcal{X}_{u,k}^v(t)>v\RT)}dv}=0.
\EQNY}
\QED

\subsection{Appendix B}
Below we state two results which are used in our proofs.

\BT\label{Cor2}
Let $\vk{X}(t), t\in[0,T]$, be a centered vector Gaussian process taking valued in $\R^n$ with continuous trajectories and $G(\vk{x}),\ \vk{x}\in\R^n$ is such that $\abs{G(\vk{x})}\leq G\abs{x}^p$ for some $G,\ p>0$ and all $\vk{x}$. Assume that the covariance matrix
\BQNY
\mathcal{R}_{t}=\E{\vk{X}(t)\vk{X}(t)^\top}
\EQNY
is uniformly non-degenerated in $[0,T]$. Assume also that for some $\Gamma, \gamma>0$,
\BQNY
\E{\abs{\vk{X}(t)-\vk{X}(s)}^2}\leq \Gamma\abs{t-s}^\gamma.
\EQNY
Denote
\BQNY
m(T)=\min_{t\in[0,T],\abs{\vk{e}}=1}\abs{\mathcal{R}_{t}^{-1/2}\vk{e}}.
\EQNY
Then there exist constants $\Gamma_1$ and $\gamma_1$ such that for all $u\geq 1$,
\BQNY
\pk{\sup_{t\in[0,T]}G(\vk{X}(t))>u}\leq \Gamma_1u^{\gamma_1}\exp\LT(-\frac{m^2(S)u^{2/\beta}}{2G^{2/\beta}}\RT).
\EQNY
\ET
This theorem is derived from \cite{PiterChaos2015}[Corollary 2].\\
 Next, we give the theorem about asymptotic expansions for probability density and tail probability of the Gaussian random chaos in
 \cite{chaos13, randomchaos15, chaos15}.\\
Denote by $g''_{d-1-m}(\vk{\varphi})$ any non-singular $(d-1-m)$-sub-matrix of $g''(\vk{\varphi})$ in \eqref{g2} and \cLb{$J(r,\vk{\varphi})$ the same as in \eqref{JJ1}}. With $\vk{\varphi}$ fixed, determinants of all such sub-matrices are the same as can be seen by applying suitable orthogonal transformations.
\BT\label{Them00}
If $g(\vk{\varphi})$ satisfies \necon{gcon} and the probability density $p_{g(\vk{\xi})}(x)$ of $g(\vk{\xi})$ exists, then the following asymptotic relations holds,
\BQN
&&p_{g(\vk{\xi})}(x)=\frac{h_0}{p\gg}\LT(\frac{x}{\gg}\RT)^{\frac{m+1}{p}-1}
\exp\LT(-\frac{x^{2/p}}{2\gg^{2/p}}\RT)(1+o(1)),\label{pdf}\\
&&\pk{g(\vk{\xi})>x}=h_0\LT(\frac{x}{\gg}\RT)^{\frac{m-1}{p}}
\exp\LT(-\frac{x^{2/p}}{2\gg^{2/p}}\RT)(1+o(1)),\label{tdf}
\EQN
as $x\rw\IF$, where for the case (i) in \necon{gcon}
\BQNY
h_0:=\frac{1}{\sqrt{2\pi}}(p\gg)^{\frac{d-1}{2}}\sum_{\vk{\varphi}\in\mathcal{M}_\varphi}
\frac{J(1,\vk{\varphi})}{\sqrt{\LT|\det g''(\vk{\varphi})\RT|}};
\EQNY
and for the case (ii) in \necon{gcon}
\BQNY
h_0:=\frac{1}{(2\pi)^{(m+1)/2}}(p\gg)^{\frac{d-1-m}{2}}\int_{\mathcal{M}_{\vk{\varphi}}}
\frac{J(1,\vk{\varphi})}{\sqrt{\LT|\det g''_{d-1-m}(\vk{\varphi})\RT|}}d V_{\vk{\varphi}},
\EQNY
with $d V_{\vk{\varphi}}$ the elementary volume in $\mathcal{M}_{\vk{\varphi}}$.
\ET
\section*{Acknowledgments} 
\COM{Partial support from the Swiss National Science Foundation Projects 200021-140633/1, 200020-159246/1 is kindly acknowledged.
K. D\c{e}bicki also acknowledges  partial support by NCN Grant No 2013/09/B/ST1/01778 (2014-2016).}
Thanks to Professor Enkelejd Hashorva and Professor Krzysztof D\c{e}bicki  who give many useful suggestions which greatly improve this manuscript. Thanks to Swiss National Science Foundation Grant no.  200021-175752.
\bibliographystyle{ieeetr}
\bibliography{Chaos}
\end{document}